\documentclass[11pt]{article}
\usepackage{amsmath,amsthm,amssymb,color,epic,eepic,
cite}
\renewcommand{\baselinestretch}{1.2}

\def\baselinestretch{1.4}
\newlength{\minitwocolumn}
\newcommand{\Z}{{\Bbb Z}} 
\newcommand{\R}{{\Bbb R}} 
\newcommand{\C}{{\Bbb C}} 
\newcommand{\N}{{\Bbb N}} 
\newcommand{\FF}{{\Bbb F}} 
\newcommand{\T}{{\Bbb T}} 

\newcommand{\F}{{\cal F}}

\newcommand{\cD}{{\cal D}}
\newcommand{\cA}{{\cal A}}

\newcommand{\cB}{{\cal B}}

\newcommand{\cI}{{\cal I}}
\newcommand{\cH}{{\cal H}}
\newcommand{\cN}{{\cal N}}
\newcommand{\cR}{{\cal R}}
\newcommand{\cP}{{\cal P}}
\newcommand{\cE}{{\cal E}}
\newcommand{\cM}{{\cal M}}
\newcommand{\cQ}{{\cal Q}}
\newcommand{\cW}{{\cal W}}

\newcommand{\bR}{{\overline{R}}}
\newcommand{\hd}{\widehat{d}}
\newcommand{\hL}{\widehat{L}}

\newcommand{\cL}{{\cal L}}

\renewcommand{\H}{{\cal H}}
\newcommand{\la}{\lambda}
\newcommand{\La}{\Lambda}

\newcommand{\al}{\alpha}

\newcommand{\ep}{\epsilon}
\newcommand{\vep}{\varepsilon}

\newcommand{\tW}{\widetilde{W}}
\newcommand{\bep}{\bar{\epsilon}}

\newcommand{\ha}{{\alpha}}

\newcommand{\hf}{\widehat{f}}

\newcommand{\hV}{\widehat{V}}
\newcommand{\hW}{\widehat{W}}

\newcommand{\bh}{{\bar{\h}}}

\newcommand{\hrho}{
{\rho}}

\newcommand{\nn}{{\nonumber}}
\newcommand{\bea}{\begin{eqnarray}}
\newcommand{\ena}{\end{eqnarray}}
\newcommand{\beit}{\begin{itemize}}
\newcommand{\enit}{\end{itemize}}

\newcommand{\be}{\begin{eqnarray*}}
\newcommand{\en}{\end{eqnarray*}}
\newcommand{\lb}[1]{\label{#1}}

\newcommand{\ds}[1]{{\displaystyle #1 }}

\newcommand{\End}{{\rm End}}

\newcommand{\id}{{\rm id}}

\def\infq4p#1{{(#1;q^4,p)_\infty}}

\newcommand{\tPhi}{\widetilde{\Phi}}

\newcommand{\tot}{\widetilde{\otimes}}

\newcommand{\tr}{{\rm tr}}

\newcommand{\mmatrix}[1]{\begin{matrix} #1 \end{matrix}}

\font\teneufm=eufm10
\font\seveneufm=eufm7
\font\fiveeufm=eufm5
\newfam\eufmfam
\textfont\eufmfam=\teneufm
\scriptfont\eufmfam=\seveneufm
\scriptscriptfont\eufmfam=\fiveeufm

\let\goth\frak
\newcommand{\slth}{\widehat{\goth{sl}}_2}
\newcommand{\slt}{\goth{sl}_2}

\newcommand{\slnh}{\widehat{\goth{sl}}_N}
\newcommand{\sln}{\goth{sl}_N}
\newcommand{\g}{\goth{g}}

\newcommand{\Bqla}{{{\cal B}_{q,\lambda}}}

\newcommand{\gl}{{\goth{gl}}}
\newcommand{\gln}{{\goth{gl}_N}}
\newcommand{\glnh}{\widehat{\goth{gl}}_N}

\newcommand{\h}{\goth{h}}
\newcommand{\gh}{\widehat{\goth{g}}}

\newcommand{\gS}{\goth{S}}

\makeatletter
\@addtoreset{equation}{section}
\makeatother

\newtheorem{thm}{Theorem}[section]
\newtheorem{prop}[thm]{Proposition}
\newtheorem{lem}[thm]{Lemma}
\newtheorem{cor}[thm]{Corollary}
\newtheorem{conj}[thm]{Conjecture}

\newtheorem{dfn}[thm]{Definition}

\begin{document}


\vspace{-1cm}
\begin{center}
{\bf\Large  Elliptic Weight Functions and  Elliptic $q$-KZ Equation
\\[7mm] }
{\large  Hitoshi Konno}\\[6mm]
{\it  Department of Mathematics, Tokyo University of Marine Science and 
Technology, \\Etchujima, Koto-ku, Tokyo 135-8533, Japan\\
       hkonno0@kaiyodai.ac.jp}
\end{center}

\begin{abstract}
\noindent 
By using representation theory of the elliptic quantum group $U_{q,p}(\slnh)$, we present a systematic method of deriving the weight functions. The resultant $\mathfrak{sl}_N$ type elliptic weight functions are new and give elliptic and dynamical analogues of those obtained 
in the rational and the trigonometric cases. We then discuss some basic properties of the elliptic weight functions. We also present an explicit formula for formal elliptic hypergeometric integral solution to the face type, i.e. dynamical, elliptic $q$-KZ equation.

\end{abstract}
\nopagebreak

\section{Introduction}
The weight functions are one of the main parts  in  hypergeometric integral solutions to the $q$-KZ equations. See for example \cite{TV97}. 
Recently new interests in the weight functions have been developing.  
Among others Gorbounov, Rim\' anyi, Tarasov and Varchenko\cite{GRTV} have succeeded to identify the rational weight functions
 with the stable envelopes 
associated with the torus equivariant cohomology of the cotangent bundle to the partial flag variety.
The stable envelopes were introduced by  Maulik and Okounkov\cite{MO,Okounkov} in more general setting associated with the 
equivariant cohomology of Nakajima's quiver variety\cite{Nakajima}. They are maps from the equivariant cohomology of the fixed point set of the torus action to the equivariant cohomology of the  variety and play an important role in a  formulation of  a geometric representation theory of quantum groups on the equivariant cohomology. 

This identification has been 
extended to the equivariant $K$-theory case in \cite{RTV} as well as   to the dynamical version of the equivariant cohomology\cite{RV} 
 and  equivariant elliptic cohomology\cite{FRV} cases both associated with the cotangent bundles to the Grassmannians. 
It has also been succeeded to construct a geometric representation 
of the Yangian $Y(\gln)$ \cite{GRTV}, the affine quantum group $U_q(\glnh)$\cite{RTV}, and  the rational dynamical quantum group $E_y(\gl_2)$\cite{RV}  and the elliptic dynamical quantum group $E_{\tau,y}(\gl_2)$\cite{FRV}  on the corresponding equivariant 
cohomology, $K$-theory,   dynamical version of the cohomology and  elliptic cohomology, respectively.  
There it  is essential to consider the finite-dimensional representations of the quantum groups on the Gelfand-Tsetlin basis of 
the tensor  product of the vector representations. 
They are lifted to the geometric representations via the correspondence between the weight functions and the stable envelopes.  
 It is also remarkable that the elliptic stable envelopes introduced by Aganagic and Okounkov\cite{AO} are the dynamical ones, where the K\"ahler variables play a role of the dynamical parameters. 

The purpose of this paper and subsequent papers is to extend these constructions
to the higher rank elliptic quantum group $U_{q,p}(\slnh)$\cite{K98,JKOS,Konno16}, which
 is an elliptic and dynamical analogue of the Drinfeld's new realization of $U_q(\slnh)$\cite{Drinfeld} and is isomorphic to 
  the central extension of Felder's elliptic quantum group\cite{Felder1,Konno16}. 
We expect that the elliptic weight functions of the $\sln$-type can be identified with the elliptic stable envelops\cite{AO}  
associated with the torus equivariant elliptic cohomology of the cotangent bundles to the partial flag variety.  
 Such identification should allow us to formulate a geometric action of 
  $U_{q,p}(\slnh)$ on the equivariant elliptic cohomology. 
Recently, Felder, Rim\' anyi and Varchenko has accomplished such study  in  the $\slth$ case\cite{FRV}. 

In this paper we discuss the elliptic weight functions of type $\sln$ and study their properties.  
We give a systematic derivation of the elliptic weight functions by using the vertex operators of  $U_{q,p}(\slnh)$\cite{JKOS, KK03, KKW}. 
For the $\slth$ case we use the level-$k$ representation of $U_{q,p}(\slth)$, and  the resultant elliptic weight functions coincide with those 
obtained in \cite{FTV1,TV97}. We discuss this case in a separate paper. We here concentrate on  the higher rank  level-1 representation  
and obtain a new result. The resultant elliptic weight functions are described by using the partitions of $[1,n]$ in a combinatorial way 
and give elliptic and dynamical analogues to those obtained by Mimachi and Noumi\cite{Mimachi,MN}, Tarasov and Varchenko \cite{TV95,RTV}. 
The same method can be applied to the trigonometric  case too by using the vertex operators of $U_q(\slnh)$ 
and allows us to derive the trigonometric weight functions  in \cite{Matsuo, TV97} for $\slth$ case and in \cite{RTV,Mimachi,MN} for 
$\sln$ case.  See also \cite{JM,TV95,Nakayashiki}. In the sequent paper we will discuss the finite-dimensional representations of $U_{q,p}(\slnh)$ on the Gelfand-Tsetlin basis and their geometric interpretations. 

Our derivation has the following  advantages. 
\begin{itemize}
\item[1)] It makes a representation theoretical meaning of the combinatorial structure of  the weight functions as well as of the partial flag variety transparent. Sec.\ref{combNot}. 
\item[2)] It makes the transition property of the weight function manifest. Sec.\ref{TransProp}.
\item[3)] It allows us to derive the shuffle algebra structure of the space of the weight functions.  Sec.\ref{shufflealgstr} $\&$ Appendix\ref{shufflealg}. 
\end{itemize}

As a byproduct we also give a new formula for formal elliptic hypergeometric integral solution to the face type elliptic $q$-KZ equation\cite{Felder2,FJMMN}.

A part of the results  has been presented in the workshops ``Classical and Quantum integrable systems'', July 11-15, 2016,  
EIMI, St.Petersburg, ``Recent Advances in Quantum Integrable Systems'', August 22-26, 2016, Univ.of Geneva and 
``Elliptic Hypergeometric Functions in Combinatorics,
Integrable Systems and Physics'', March 20-24, 2017, ESI, Vienna.

This paper is organized as follows. 
In Section 2 we prepare some notations including the elliptic dynamical $R$ matrices. 
In Section 3 we review a construction of the vertex operators of the $U_{q,p}(\slnh)$-modules. 
In particular we provide a free field realization of the vertex operators for the  level 1 representation. 
Section 4 is devoted to a derivation of the elliptic weight functions of the $\sln$-type by using the vertex operators. 
In Section 5 we discuss some basic  properties of the elliptic weight functions such as the triangular property,  transition property, 
orthogonality, quasi-periodicity and the shuffle algebra structures.  
 In Section 6, we give a formal elliptic hypergeometric integral solution to the face type elliptic $q$-KZ equation. 
In Appendix A we summarize some basic facts on the elliptic quantum group $U_{q,p}(\slnh)$. 
Appendix B is a summary on the dynamical representations of $U_{q,p}(\slnh)$ both the finite and infinite-dimensional cases.  
In Appendix C we provide a proof of Proposition \ref{shuffleMp} including a derivation of the shuffle algebra structure of the space of 
the weight functions. 

While preparing this manuscript, we become aware of the paper by Rim\' anyi, V.Tarasov and A.Varchenko \cite{RTV17} which has 
some overlap with the content of the present paper. 

\section{Preliminaries}
In this section we prepare the notation to be used in the text. 
Throughout this paper  $q$ are  generic complex numbers satisfying $|q|<1$ unless otherwise specified.

\subsection{The commutative algebra $H$}
Let  $A=(a_{ij})\ (0\leq i,j\leq N-1)$ be   
 the  generalized Cartan matrix of $\slnh=\widehat{\mathfrak{ sl}}(N,\C)$\cite{Kac}. 
Let $\h=\widetilde{\h}\oplus \C d$, $\widetilde{\h}=\bar{\h}\oplus\C c$, $\bar{\h}=\oplus_{i=1}^{N-1}\C h_i$ be the Cartan subalgebra of $\slnh$.  
 Define  $\delta, \Lambda_0, \al_i, \bar{\Lambda}_i\ (1\leq i\leq N-1) \in \h^*$ 
by 
\bea
&&<\al_i,h_j>=a_{ji}, \ <\delta,d>=1=<\Lambda_0,c>,\ <\bar{\Lambda}_i,h_j>=\delta_{i,j}\lb{pairinghhs}
\ena
the other pairings are 0. 
We set $\bh^*=\oplus_{i=1}^{N-1} \C \bar{\Lambda}_i, $ $\widetilde{\h}^*=\bh^*\oplus \C \Lambda_0$, 
 $\cQ=\oplus_{i=1}^{N-1}\Z \al_i$ and 
$\cP=\oplus_{i=1}^{N-1}\Z \bar{\Lambda}_i$. 
 Let $\{ \ep_j\ (1\leq j\leq N)\}$ be an orthonormal basis  in $\R^N$ 
with the inner product
$( \ep_j, \ep_k )=\delta_{j,k}$.  We set $\bep_j=\ep_j-\sum_{k=1}^N\ep_k/N$. 
We realize the simple roots by $\alpha_j=\ep_j-\ep_{j+1}\ (1\leq j\leq N-1)$ and the fundamental weights by $\bar{\Lambda}_j=\bep_1+\cdots+\bep_j \ (1\leq j\leq N-1)$. 
We define $h_{\ep_j}\in \bar{\h} \ (1\leq j\leq N)$ by 
$<\ep_i,h_{\ep_j}>=(\ep_i,{\ep_j})$ and $h_\alpha \in \bar{\h}$ for $\alpha=\sum_j c_j \ep_j, \ (c_j\in \C)$ by
 $h_\alpha=\sum_j c_j h_{\ep_j}$.
We regard $\bar{\h}\oplus \bar{\h}^*$ as the Heisenberg 
algebra  by
\bea
&&~[h_{{\epsilon}_j},{\epsilon}_k]
=( {\epsilon}_j,{\epsilon}_k ),\qquad [h_{{\epsilon}_j},h_{{\epsilon}_k}]=0=
[{\epsilon}_j,{\epsilon}_k].\lb{HA1}
\ena
In particular, we have $[h_{j}, \alpha_k]=a_{j k}$.  We also set $h^j=h_{\bar{\Lambda}_j}$. 
 
Let $\{P_{{\alpha}}, Q_{{\beta}}\}\ 
({\alpha}, {\beta} \in \bar{\h}^*)$ be the Heisenberg algebra 
defined by the commutation relations
\begin{eqnarray}
&&[P_{{\epsilon}_j}, Q_{{\epsilon}_k}]=
( {\epsilon}_j, {\epsilon}_k ), \qquad 
[P_{{\epsilon}_j}, P_{{\epsilon}_k}]=0=
[Q_{{\epsilon}_j}, Q_{{\epsilon}_k}],\lb{HA2}
\end{eqnarray}
where
$P_\alpha=\sum_j c_j P_{\ep_j}$ for $\alpha=\sum_j c_j \ep_j$.
We set $P_{\bh}=\oplus_{j=1}^{N}\C P_{\ep_j}, Q_{\bh}=\oplus_{j=1}^{N}\C Q_{\ep_j}$  $P_{j}=P_{\alpha_j}, P^j=P_{\bar{\Lambda}_j}$ and 
$Q_{j}=Q_{\alpha_j}, Q^j=Q_{\bar{\Lambda}_j^\vee}$. 

For the abelian group 
$\cR_Q= \oplus_{i=1}^{N-1}\Z Q_{\al_i}$, 
we denote by $\C[\cR_Q]$ the group algebra over $\C$ of $\cR_Q$. 
We denote by $e^{Q_\al}$ the element of $\C[\cR_Q]$ corresponding to $Q_\al\in \cR_Q$. 

We define the  commutative algebra $H$ by $H=\widetilde{\h}\oplus P_{\bh}=\sum_{j}\C(P_{\ep_j}+h_{\ep_j})+\sum_j\C P_{\ep_j}+\C c$. We  denote the dual space of $H$ by $H^*=\widetilde{\h}^*\oplus Q_{\bh}$. We define the paring by \eqref{pairinghhs},  
$<Q_\al,P_\beta>=(\al,\beta)$ and  $<Q_\al,h_\beta>=<Q_\al,c>=<Q_\al,d>=0=<\al,P_\beta>=<\delta,P_\beta>=<\Lambda_0,P_\beta>$. 
We often abbreviate $P_{\ep_j}+h_{\ep_j}$ as  $(P+h)_{\ep_j}$ and use $(P+h)_{j,k}=(P+h)_{\ep_j}-(P+h)_{\ep_k}$, $P_{j,k}=P_{\ep_j}-P_{\ep_k}$, $h_{j,k}=h_{\ep_j}-h_{\ep_k}$ etc..  

We define $\FF=\cM_{H^*}$ to be the field of meromorphic functions on $H^*$. 

\subsection{Infinite products}
We use the following notations.
\be 
&&[n]_q=\frac{q^n-q^{-n}}{q-q^{-1}},
 \qquad  (x;q)_\infty=\prod_{n=0}^\infty(1-x q^n),\\
&& (x;q,t)_\infty=\prod_{m,n=0}^\infty(1-x q^m t^n),\quad  (x;p,q,t)_\infty=\prod_{l,m,n=0}^\infty(1-x p^lq^m t^n),\\
&&\Theta_p(z)=(z;p)_{\infty}(p/z;p)_\infty(p;p)_\infty, \\
&&\Gamma(x;q,t)=\frac{(qt/x;q,t)_\infty}{(x;q,t)_{\infty}}, \qquad \Gamma(x_1,x_2;q,t)=\Gamma(x_1;q,t)\Gamma(x_2;q,t),\\
&&\Gamma(x;p,q,t)=(x;p,q,t)_\infty(pqt/x;p,q,t)_\infty.
\en
for $|q|<1, |t|<1, |p|<1$.

\subsection{Theta functions}
Let $r$ be a generic positive real number. 
We set $p=q^{2r}$.  In general we consider the level $k\in \R$ representation of 
the elliptic quantum group $U_{q,p}(\g)$.  
See Appendix \ref{App:rep}.  In that case we assume $r^*=r-k>0$ and set $ p^*=q^{2r^*}$. 

 We use the following Jacobi's odd theta functions.
\bea
&&[u]=q^{\frac{u^2}{r}-u}\Theta_p(z), \qquad  
[u]^*=q^{\frac{u^2}{r^*}-u}\Theta_{p^*}(z), \lb{thetafull}\\
&&[u+r]=-[u], \quad [u+r\tau]=-e^{-\pi i\tau}e^{-2\pi i {u}/{r}}[u],  \quad \lb{thetaquasiperiod}\\
&&[u+r^*]^*=-[u]^*, \quad [u+r^*\tau^*]^*=-e^{-\pi i\tau^*}e^{-2\pi i {u}/{r^*}}[u]^*, \nn
\ena
where $z=q^{2u}, p=e^{-2\pi i/\tau }, p^*=e^{-2\pi i/\tau^* }$.

\subsection{The elliptic dynamical $R$-matrix of  type $\slnh$}\lb{R-mat}

Let $\hV_z$ be the $N$-dimensional dynamical  evaluation representation of $U_{q,p}(\slnh)$ given in Sec \ref{vecrep} 
and $\{v_{\mu}\ (\mu=1,\cdots,N)\}$ be its basis. 
We consider the following elliptic dynamical $R$-matrix  $R^+(z_1/z_2,\Pi)\in \End_{\C}(\hV_{z_1}\tot \hV_{z_2})$ given by
\begin{eqnarray}
R^+(z,\Pi)&=&\hrho^+(z){\bR}(z,\Pi),\lb{def:Rmat}\\
{\bR}(z,\Pi)&=&
\sum_{j=1}^{N}E_{j,j}\otimes E_{j,j}+
\sum_{1 \leq j_1< j_2 \leq N}
\left(b_{}(u,(P+h)_{j_1,j_2 })
E_{j_1,j_1}
\otimes E_{j_2,j_2}+
\bar{b}_{}(z)
E_{j_2,j_2}\otimes E_{j_1,j_1}
\right.
 \nonumber\\
&&\qquad 
\left.
+
c_{}(u,(P+h)_{j_1,j_2 })
E_{j_1,j_2}\otimes E_{j_2,j_1}+
\bar{c}_{}
(u,(P+h)_{j_1,j_2 })E_{j_2,j_1}\otimes E_{j_1,j_2}
\right),\lb{ellR}
\end{eqnarray}
where $z=q^{2u}$, $\Pi_{j,k}=q^{2(P+h)_{j,k}}$, 
\bea
&&{\rho}^+(z)=q^{-\frac{N-1}{N}}z^{\frac{N-1}{rN}}\tilde{\rho}^{+}(z),\qquad 
\tilde{\rho}^+(z)=\frac{\{q^{2N}q^{-2} z\}\{q^2 z\}}{\{q^{2N} z\}\{z\}}
\frac{\{pq^{2N}/z\}\{p/z\}}{\{pq^{2N}q^{-2}/z\}\{pq^2/z\}},\lb{defrhotilde}\\ 
&&b(u,s)=
\frac{[s+1][s-1][u]}{[s]^2[u+1]},\qquad 
\bar{b}(u)=
\frac{[u]}{[u+1]},\lb{Relements}\\
&&c(u,s)=\frac{[1][s+u]}{[s][u+1]},\qquad \bar{c}(u,s)=\frac{[1][s-u]}{[s][u+1]}\nn
\ena
and $\{z\}=(z;p,q^{2N})_\infty$. 
This $R$ matrix is gauge equivalent to Jimbo-Miwa-Okado's $A_{N-1}^{(1)}$ face type Boltzmann weight\cite{JMO} 
and can be obtained by taking the vector representation of the universal elliptic dynamical $R$ matrix derived in \cite{JKOStg}.
See \cite{Konno06}. 

We also use $R^{*+}(z,\Pi^*)$ with $\Pi^*_{j,k}=q^{2P_{j,k}}$, which is obtained from $R^{+}(z,\Pi)$ by the 
following replacement.   
\be
&&R^{*+}(z,\Pi^*)=\left.R^{+}(z,\Pi)\right|_{p\to p^*, r \to r^*, [u]\to [u]^*, P+h \to P}.
\en

The $R^+(z,q^{2s})$ satisfies the dynamical Yang-Baxter equation
\bea
&& R^{+(12)}(z_1/z_2,q^{2(s+h^{(3)})}) R^{+(13)}(z_1/z_3,q^{2s}) R^{+(23)}(z_2/z_3,q^{2(s+h^{(1))})})\nn\\
&&\qquad =R^{+(23)}(z_2/z_3,q^{2s}) R^{+(13)}(z_1/z_3,q^{2(s+h^{(2)})}) R^{+(12)}(z_1/z_2,q^{2s}).\lb{DYBE}
\ena
where  $q^{2h_{j,k}^{(l)}}$ acts on the $l$-th tensor space $\hV_{z_l}$ by  $q^{2h_{j,k}^{(l)}}  v_{\mu}=q^{2<\bep_\mu,h_{j,k}>}v_\mu$,　
and the unitarity
\bea
&&R^+(z,q^{2s})R^{(21)}(z^{-1},q^{2s})=\id_{\hV_z\otimes \hV_1}.
\ena

\section{Vertex Operators of $U_{q,p}(\slnh)$}\lb{VO}
In this section we summarize the known facts on  the type I and II vertex operators of the $U_{q,p}(\slnh)$-modules 
obtained in \cite{JKOS,KK03}.

\subsection{Definition}
Let $\widehat{V}_z$ be as in Sec.\ref{R-mat} 
and $\widehat{V}(\lambda,\nu)$ denote the irreducible 
level-$k$ highest weight $U_{q,p}(\slnh)$-module with highest weight $(\lambda,\nu)$ in Definition \ref{def:levelkrep}. 
The level-1 case is given in Theorem \ref{levelone}. 
The type $\mathrm{I}$ $\Phi(z)$ and the type $\mathrm{II}$ $\Psi^*(z)$ vertex operators are the intertwiners of 
the $U_{q,p}(\slnh)$-modules of the form
\bea
 {\Phi}(z) &:& \widehat{V}(\lambda,\nu) \to  \widehat{V}_z \tot \widehat{V}(\lambda',\nu), \\
  {\Psi}^*(z) &:&  \widehat{V}(\lambda,\nu) \tot \widehat{V}_z\to \widehat{V}(\lambda,\nu'),
\ena
where $\la, \la'\in \h^*,  \nu, \nu' \in H^*$. 
The vertex operators satisfy the intertwining relations with respect to the comultiplication
$\Delta$ given in Sec.\ref{Hopfalgebroid}
\bea
 && \Delta(x) {\Phi}(z) ={\Phi}(z)x,         \label{irI}\\
 && x {\Psi}^*(z) = {\Psi}^*(z) \Delta(x),   \qquad \forall x \in U_{q,p}(\slnh).\label{irII}
\ena
By using the $L$ operator $\hL^+(z)$ given in \eqref{def:lhat} the main part of the intertwining relations can be re-expressed 
as  follows\cite{Konno08}:
\bea
&&(\id\tot {\Phi}(z_2))\widehat{L}^{+}(z_1)=
   R^{+(12)}(z_1/z_2,\Pi)\widehat{L}^{+(13)}(z_1)(\id \tot {\Phi}(z_2)),
 \label{typeI}
\\
&& \widehat{L}^{+}(z_1) {\Psi}^{*(23)}(z_2)=
   {\Psi}^{*(23)}(z_2)\widehat{L}^{+(12)}(z_1)R^{+*(13)}(z_1/z_2,\Pi^*q^{-2(h^{(1)}+h^{(3))}}).
\label{typeII}
\ena
The relation (\ref{typeI}) (resp.  (\ref{typeII})) should be understood on $\widehat{V}_{z_1} \tot \widehat{V}(\lambda,\nu)$ (resp. $\widehat{V}_{z_1} \tot \widehat{V}(\lambda,\nu) \tot  \widehat{V}_{z_2}$). 

We define the components of the vertex operators by
\begin{equation}
 {\Phi}(zq^{-1})u=\sum_{\mu=1}^{ N}v_{\mu} 
 \tot \Phi_{\mu} \left(z \right)u, \quad
 {\Psi}^*(zq^{k-1})( u\tot v_{\mu})=\Psi_{\mu}^* \left(z \right) u \qquad \forall u \in \widehat{V}(\la,\nu).
\end{equation}

Let $F_j(z), E_j(z), K^+_j(z)$ be the elliptic currents in Definition \ref{modellcurrents} and let $F^+_{j,l}(z), E^+_{l,j}(z)$ be the 
half currents in Definition \ref{def:halfcurrents}. 
From the intertwining relations \eqref{typeI}-\eqref{typeII} one can deduce the following relations as the sufficient conditions\cite{KK03}. 
\begin{prop}\lb{sufvo}
For the type I, 
\bea
&& \Phi_\mu(z_2)=F_{\mu,N}^+(pq^{-1}z_2)\Phi_{N}(z_2)\qquad (1 \leq \mu \leq N-1),\lb{typeIk}
\ena
and
\bea
&& \Phi_{N}(z_2)K_{N}^+(z_1)=\rho^+(qz_1/z_2)K_{N}^+(z_1)\Phi_{N}(z_2),\\
&& [\Phi_{N}(z),P_{\bep_l}]=0, \quad [\Phi_{N}(z),E_j(w)]=0. \qquad (1\leq l\leq N, 1\leq j\leq N-1),\\
&& [\Phi_{N}(z),(P+h)_{k,l}]=-\delta_{l,N}\Phi_{N}(z),\quad (k<l)
 \label{typIsufcnd3}\\
&& F_{N-1}(z_1)\Phi_{N}(z_2)= \frac{[u_1-u_2+\frac{1}{2}]}
  {[u_1-u_2-\frac{1}{2}]}\Phi_{N}(z_2)F_{N-1}(z_1),\label{PhiNFNm1}\\
&&F_j(z_1)\Phi_{N}(z_2)=\Phi_{N}(z_2)F_j(z_1) \quad (1 \le j \le N-2).
\label{typIsufcnd5}
\ena
For the type II, 
\bea
&&\Psi_\mu^*(z_2)=\Psi^*_{N}(z_2)E_{N,\mu}^+(p^*z_2)
\quad (1 \leq \mu \leq N-1),
\ena
and
\bea
&& K_{N}^+(z_1)\Psi_{N}^*(z_2)=\Psi_{N}^*(z_2)K_{N}^+(z_1)\rho^{+*}(q^{1-k}z_1/z_2),\label{typIIsufcnd1}\\
&& [ \Psi_{N}^*(z_2),(P+h)_{\bep_l}]=0, \quad [\Psi_{N}^*(z_2),F_j(z_1)]=0\qquad (1\leq l\leq N, 1\leq j\leq N-1), 
 \label{typIIsufcnd2}\\
&& [\Psi_{N}^*(z),P_{k,l}]=-\delta_{l,N}\Psi_{N}^*(z)\qquad  (k<l),
 \label{typIIsufcnd3}\\
&& E_{N-1}(z_1)\Psi_{N}^*(z_2)=\frac{[u_1-u_2-\frac{1}{2}]^*}{[u_1-u_2+\frac{1}{2}]^*}
 \Psi_{N}^*(z_2)E_{N-1}(z_1),
 \label{typIIsufcnd4}\\
&&E_j(z_1)\Psi_{N}^*(z_2)=\Psi_{N}^*(z_2)E_j(z_1) \quad (1 \le j \le N-2).
\label{typIIsufcnd5}
\ena

\end{prop}

\subsection{The level-1 vertex operators and commutation relations}
For the level-1 representation (Sec.B.3),  we have a free field realization of the vertex operators.

\begin{thm}\lb{bareboson}\cite{KK03}
Let us assume $|p|<|z|<1$. Let $\Lambda_a\ (a=0,1,\cdots,N-1)$ be the fundamental weights of $\slnh$.  The components of the type I  and the type II  vertex operators,
$\Phi_\mu(z): \hV(\La_a+\nu,\nu)\to \hV(\La_a-\mu+\nu,\nu)$ and 
$\Psi^*_\mu(z): \hV(\La_a+\nu,\nu)\to \hV(\La_a+\nu,\nu-\mu)$,  are realized as follows. 
\bea
 \Phi_{N}(z) &=& : \exp \left\{ \sum_{m \ne 0} (q^m-q^{-m}){\cE'}_m^{+N} z^{-m} \right\}:
  e^{-\widehat{\bep}_N
} (- z)^{-h_{\bep_N}}z^{\frac{1}{r}(P+h)_{\bep_N}}
 , \label{typeIvo}
 \\
\Phi_{\mu}(z)&=&F^+_{\mu,N}(q^{-1}pz)\Phi_{N}(z)\nn\\
&=&a_{\mu,N}\oint_{\T^{N-\mu}}\prod_{m=\mu}^{N-1}\frac{dt_m}{2\pi i t_m}\Phi_{N}(z)F_{N-1}(t_{N-1})F_{N-2}(t_{N-2})\cdots F_{\mu}(t_{\mu})
\varphi_{\mu}(z,t_{\mu},\cdots,t_{N-1};\Pi),\nn\\
&&\hspace{5cm}(1\leq \mu\leq N-1),\lb{typeImu}
\\
 \Psi_{N}^*(z) &=&: \exp \left\{ -\sum_{m \ne 0} (q^m-q^{-m})
   \cE_m^{+N}z^{-m} \right\}: 
  e^{\widehat{\bep}_N
}e^{-Q_{\bep_N}}(- z)^{h_{\bep_N}} z^{-\frac{1}{r^*}P_{\bep_N}},
   \label{typeIIvo}
\ena
\bea
 \Psi_{\mu}^*(z)&=& \Psi_{N}^*(z)E^+_{N,\mu}(p^*z)\nn\\
 &=&a^*_{\mu,N}\oint_{\T^{N-\mu}}\prod_{m=\mu}^{N-1}\frac{dt_m}{2\pi i t_m}\Psi^*_{N}(z)E_{N-1}(t_{N-1})E_{N-2}(t_{N-2})\cdots E_{\mu}(t_{\mu})
\varphi^*_{\mu}(z,t_{\mu},\cdots,t_{N-1};\Pi^*),\nn\\
&&\hspace{5cm}(1\leq \mu\leq N-1)
\ena
where  
$\cE_m^{+j}$ and $\displaystyle {\cE'}_m^{+j}=\frac{1-p^{*m}}{1-p^m}q^{m}\cE_m^{+j}$ are 
the orthonormal basis type elliptic bosons of level 1  
 given in \eqref{orthoboson}.  
We also set $\widehat{\bep}_N=\bep_N+\zeta_{\bep_N}$ as  in Sec.\ref{levelonerep}, 
 $\Pi=\{\Pi_{\mu,m}\ (m=\mu+1,\cdots,N)\}
$, $\Pi^*=\{\Pi^*_{\mu,m}\ (m=\mu+1,\cdots,N)\}
$, 
\bea
&&{\varphi_{\mu}(z,t_{\mu},\cdots,t_{N-1}; \Pi)=
\prod_{m=\mu}^{N-1}
\frac{[v_{m+1}-v_{m}+(P+h)_{\mu,m+1}-\frac{1}{2}][1]}
{[v_{m+1}-v_{m}+\frac{1}{2}][(P+h)_{\mu,m+1}]}}, \lb{defvarphi}\\
&&{\varphi^*_{\mu}(z,t_{\mu},\cdots,t_{N-1}; \Pi^*)=
\prod_{m=\mu}^{N-1}
\frac{[v_{m+1}-v_{m}-P_{\mu,m+1}+\frac{1}{2}]^*[1]^*}
{[v_{m+1}-v_{m}-\frac{1}{2}]^*[P_{\mu,m+1}-1]^*}} \lb{defvarphistar}
\ena
with $z=q^{2u}, t_m=q^{2v_m}$, $v_N=u$ and
\be
&&
\T^{N-\mu}=\{t\in \C^{N-\mu}\ |\ |t_\mu|=\cdots=|t_{N-1}|=1 \}.   
\en
   
\end{thm}

\begin{thm}\cite{KK03}
The free field realizations of the type I $\Phi_\mu(z)$ and the type II  $\Psi_\mu^*(z)$ vertex operators satisfy the following commutation relations:
\begin{eqnarray}
 \Phi_{\mu_2}(z_2)\Phi_{\mu_1}(z_1) &=& \sum_{\mu_1',\mu_2'=1}^{N}R(z_1/z_2,\Pi)_{\mu_1 \mu_2}^{\mu_1' \mu_2'}\ 
 \Phi_{\mu_1'}(z_1)\Phi_{\mu_2'}(z_2), \label{typeIcr}\\
  \Psi_{\mu_1}^*(z_1)\Psi_{\mu_2}^*(u_2) &=& \sum_{\mu_1',\mu_2'=1}^{N}\Psi_{\mu_2'}^*(z_2)\Psi_{\mu_1'}^*(z_1)
  R^{*}(z_1/z_2,\Pi^*)_{\mu_1' \mu_2'}^{\mu_1 \mu_2}, \label{typeIIcr} \\
\Phi_\mu(z_1)\Psi_\nu^*(z_2) &=& \chi(z_1/z_2)\Psi_\nu^*(z_2)\Phi_\mu(z_1). \label{typeI-II} 
\end{eqnarray}
Here we set
\begin{equation}
 R(z,\Pi)=\mu(z)\bR(z,\Pi), \quad R^*(z,\Pi^*)=\mu^*(z)\bR^{*}(z,\Pi^*)
\end{equation}
with
\bea
&& \mu(z)=z^{-\frac{r-1}{r}\frac{N-1}{N}} \frac{ \{p q^{2N} q^{-2} z\} 
 \{q^2 z \}}{ \{p z \} \{q^{2N} z \} } 
  \frac{  \{p /z \} 
  \{q^{2N} /z \}  }{\{p q^{2N} q^{-2}/z \} \{q^2 /z \} }, \lb{def:mu}\\
&&   \mu^*(z)=z^{\frac{r}{r-1}\frac{N-1}{N}} \frac{ \{ q^{2N} q^{-2} z\}^* 
 \{p^* q^2 z \}^*}{ \{ z \}^* \{p^*q^{2N} z \}^* } 
  \frac{  \{1 /z \}^*  \{p^*q^{2N} /z \}^*  }{\{q^{2N} q^{-2}/z \}^* \{p^*q^2 /z \}^* }\lb{def:mus}
\ena
and
\bea
&&\chi(z)=z^{-\frac{N-1}{N}}\frac{\Theta_{q^{2N}}(qz)}{\Theta_{q^{2N}}(q/z)}.
\ena
\end{thm}

\section{Elliptic Weight Functions}\lb{EWF}
In this section we present a simple derivation of the elliptic weight functions of  type $\sln$ 
 by using the vertex operators prepared in the last section.
The method can also be applied to the rational and trigonometric cases if one has an appropriate realization of the vertex operators.
 
Let us consider the $U_{q,p}(\slnh)$ level-1 type I vertex operators $\Phi_\mu(z)\ (\mu=1,\cdots,N)$ in Theorem \ref{bareboson} 
and their $n$-point composition  $\phi(z_1,\cdots,z_n)=\Phi(z_1)\circ \cdots\circ\Phi(z_n)\ : \hV(\La_a+\nu,\nu) \to \hV_{z_n}\tot\cdots\tot \hV_{z_1} $ $\tot
\hV(\La_{a'}+\nu,\nu)$, where $a'=(0 N-1 N-2 \cdots 2 1)^n(a)$ is the cyclic permutation of $a$. We have
\be
&&\phi(z_1,\cdots,z_n)=\sum_{\mu_1,\cdots,\mu_n=1}^Nv_{\mu_n}\tot\cdots\tot v_{\mu_1}\tot \phi_{\mu_1\cdots\mu_n}(z_1,\cdots,z_n),\\
&&\phi_{\mu_1\cdots\mu_n}(z_1,\cdots,z_n)=\Phi_{\mu_1}(z_1)\cdots\Phi_{\mu_n}(z_n).
\en

\subsection{Combinatorial notations}\lb{combNot}
For a given $n$-point operator $\phi_{\mu_1\cdots\mu_n}(z_1,\cdots,z_n)$, it is convenient to introduce
 the following combinatorial notations. Let $[1,n]=\{1,\cdots,n\}$. 
Define the index set $I_l:=\{\ i\in [1,n]\ |\ \mu_i=l\}$ $(l=1,\cdots,N)$ and set  $\la_l:=|I_l|$, 
$\la:=(\la_1,\cdots,\la_N)$.  
 Then $I=(I_1,\cdots,I_N)$ is a partition of $[1,n]$, i.e.  
 \be
I_1\cup \cdots \cup I_N=[1,n],\quad I_k\cap I_l=\emptyset\quad  \mbox{$(k\not=l)$}.
\en
We often denote thus obtained partition $I$ by $I_{\mu_1,\cdots\mu_n}$ and conversely the $n$-point operator by $\phi_I(z_1,\cdots,z_n)$. 
Let $\N=\{m\in \Z|\ m\geq 0 \}$. 
For $\la=(\la_1,\cdots,\la_N)\in \N^N$, let $\cI_\la$ be the set of all partitions  $I=(I_1,\cdots,I_N)$ satisfying $|I_l|=\la_l\ (l=1,\cdots,N)$.
Note that for all $I\in \cI_\la$ the $n$-point operators $\phi_{I}(z_1,\cdots,z_n)$ have the same weight $-\sum_{j=1}^n\bep_{\mu_j}$. 
We call $\sum_{j=1}^n\bep_{\mu_j}$ the weight associated with $\la$. 
 We also set $\la^{(l)}:=\la_1+\cdots+\la_l$,  $I^{(l)}:=I_1\cup\cdots \cup I_l$ and let $I^{(l)}=:\{i^{(l)}_1< \cdots<i^{(l)}_{\la^{(l)}}\}$. 
 It is also important to specify the numbering of the arguments of the elliptic currents $F_l$'s appearing in the realization of the vertex 
 operators $\Phi_\mu(z)$.  For $l=1,\cdots,N-1$, we assign the argument $t^{(l)}_k$  to  the elliptic current $F_l$ attached to the $i^{(l)}_k$-th vertex operator.

\noindent
{\it Example 1.}\ Let us consider the $N=3$, $n=5$,  $\la=(2,2,1)$ case. For example, the 5-point operator
\bea
&&\Phi_2(z_1)\Phi_1(z_2)\Phi_3(z_3)\Phi_1(z_4)\Phi_2(z_5)\lb{exVO}
\ena
gives a partition $I=(I_1=\{2,4\}, I_2=\{1,5\}, I_3=\{3\})$.  Hence $I^{(1)}=\{2,4\}$, $I^{(2)}=\{1,2,4,5\}$, $I^{(3)}=\{1,2,3,4,5\}$. 
In particular, $i^{(2)}_1=i^{(3)}_1=1$, $i^{(1)}_1=i^{(2)}_2=i^{(3)}_2=2$, $i^{(3)}_3=3$,  $i^{(1)}_2=i^{(2)}_3=i^{(3)}_4=4$,  
$i^{(2)}_4=i^{(3)}_5=5$.  Then from Theorem \ref{bareboson} we obtain the following  realization of  the vertex operators in  \eqref{exVO}.
\be
&&\Phi_2(z_1)=a_{2,3}\oint_{\T} \frac{dt^{(2)}_1}{2\pi i t^{(2)}_1}\Phi_3(z_1)F_2(t^{(2)}_1)\varphi_2(z_1,t^{(2)}_1;\Pi),\\
&&\Phi_1(z_2)=a_{1,3}\oint_{\T^2} \frac{dt^{(2)}_2}{2\pi i t^{(2)}_2}\frac{dt^{(1)}_1}{2\pi i t^{(1)}_1}\Phi_3(z_2)F_2(t^{(2)}_2)F_1(t^{(1)}_1)
\varphi_1(z_2,t^{(2)}_2,t^{(1)}_1;\Pi),\\
&&\Phi_1(z_4)=a_{1,3}\oint_{\T^2} \frac{dt^{(2)}_3}{2\pi i t^{(2)}_3}\frac{dt^{(1)}_2}{2\pi i t^{(1)}_2}\Phi_3(z_4)F_2(t^{(2)}_3)F_1(t^{(1)}_2)
\varphi_1(z_2, t^{(2)}_3,t^{(1)}_2;\Pi),\\
&&\Phi_2(z_5)=a_{2,3}\oint_{\T} \frac{dt^{(2)}_4}{2\pi i t^{(2)}_4}\Phi_3(z_5)F_2(t^{(2)}_4)\varphi_2(z_4,t^{(2)}_4;\Pi).
\en

\noindent
{\it Remark.}\cite{GRTV,RTV}\ For a partition $I$ and associated assignment of the variables $t^{(l)}_a$,  it is convenient to make a $n\times N$ table  
of the variables $z_j$ and $t^{(l)}_a$ $(j=1,\cdots,n, l=1,\cdots,N-1, a=1,\cdots,,\la^{(l)} )$ obtained  by the following rule: put $t^{(l)}_a$ into the $(i^{(l)}_a,l)$ box (i.e. a box located at the $i^{(l)}_a$-th row and the $l$-th column ) and  put $z_1,\cdots, z_n$ into the $N$-th column from top to bottom. 

\noindent 
{\it Example 2.}   The case given in {\it Example 1}, we obtain the following table.  
$$
\begin{array}{|c|c|c|}\hline
              &t^{(2)}_1&z_1\\ \hline
t^{(1)}_1&t^{(2)}_2&z_2\\ \hline
              &             &z_3\\ \hline
t^{(1)}_2&t^{(2)}_3&z_4\\ \hline
           &t^{(2)}_4&z_5\\ \hline  
\end{array}       
$$

\noindent
{\it Remark.}\cite{MN,GRTV,RTV}\ 
For $\la=(\la_1,\cdots,\la_N)\in \N^N$, the above notation is well fit for a parametrization of the partial flag variety  $\F_\la$  consisting of 
$0=\F_0\subset \F_1\subset \cdots\subset \F_N=\C^n$
with $\dim \F_l/\F_{l-1}=\la_l$.  Our assignment hence gives a representation theoretical meaning to such parametrization.

\subsection{Derivation}
Substituting the expressions of the vertex operators \eqref{typeImu} into the $n$-point operator, we obtain 
\bea
&&\phi_{\mu_1\cdots\mu_n}(z_1,\cdots,z_n)
\nn\\
&&=\oint_{\T^{N-\mu_1}}  \underline{dt_1}\  \Phi_N(z_1)F_{N-1}(t^{(N-1)}_{1})\cdots F_{\mu_1}(t^{(\mu_1)}_1){\varphi_{\mu_1}(z_1,t^{(N-1)}_{1},\cdots,t^{(\mu_1)}_1;\{\Pi_{\mu_1,m}\})}\cdots\nn \\
&&\qquad\qquad \cdots  \oint_{\T^{N-\mu_n}}  \underline{dt_n}\ \Phi_N(z_n)F_{N-1}(t^{(N-1)}_{\la^{(N-1)}})\cdots F_{\mu_n}(t^{(\mu_n)}_{\la^{(\mu_n)}}){\varphi_{\mu_n}(z_1,t^{(N-1)}_{\la^{(N-1)}},\cdots,t^{(\mu_n)}_{\la^{(\mu_n)}};\{\Pi_{\mu_n,m}\})},\nn\\
&&\lb{nptoperator}
\ena
where we set
\be
&&\underline{dt_1}=a_{\mu_1,N}
\prod_{j=\mu_1}^{N-1}\frac{d t^{(j)}_{1}}{2\pi i  t^{(j)}_{1}},\qquad 
\underline{dt_n}=a_{\mu_n,N}
\prod_{j=\mu_n}^{N-1}\frac{d t^{(j)}_{\la^{(j)}}}{2\pi i  t^{(j)}_{\la^{(j)}}}, \qquad \mbox{etc.}
\en

We then divide the integrand into two parts, 
the operator part $\tPhi(t,z)$ and the kinematical factor part $\omega_{\mu_1\cdots\mu_n}(t,z,\Pi)$.  
The operator part consists of the normal ordered bare vertex operators $\Phi_N(z)$'s,  the index-wise normal ordered elliptic currents $F_l(t)\ (l=1,\cdots,N-1)$ and the symmetric part of the OPE coefficients. 
We put  $\Phi_N(z)$'s and  $F_l(t)$'s $(l=1,\cdots,N-1)$ in the definite ordering  specified below.   
The kinematical factor part consists of all $\varphi_{\mu}(z,t_{\mu},\cdots;\Pi)$'s,  all factors arising from the exchange relations between  
$\Phi_N(z)$'s and $F_{N-1}(t)$'s  as well as  among $F_{l}(t)$'s $(l=1,\cdots,N-1)$ and of all the non-symmetric part of the OPE coefficients. 

The following procedure yields the weight functions, which are natural elliptic and dynamical analogues of the trigonometric ones obtained in \cite{Mimachi,RTV,MN}, as the kinematical factor part.
The procedure consists of the following 4 steps. 
\begin{itemize}
\item[1.] Move all $\varphi_{\mu}$'s to the right end. Then the dynamical parameters in $\varphi_\mu$ get shift 
following the exchange relation
\be
&&\Pi_{\nu,m}\Phi_{\mu}(z)=\Phi_{\mu}(z)\Pi_{\nu,m}q^{-2<\bep_{\mu},h_{\nu,m}>}. 
\en 
\item[2.] Move all the elliptic currents $F_l(t)$'s to the right of all the bare vertex operators $\Phi_N(z)$'s and arrange the order of all 
$F_l(t)$'s by collecting them into $N-1$ groups by their indices.   Put these $N-1$ groups in the decreasing order. 
Then one gets appropriate  factors by the exchange relations \eqref{u8} and \eqref{PhiNFNm1}.  
\item[3.] Take normal ordering of all $\Phi_N(z)$'s and each group  of $F_l(t)$'s having the same indices $(l=1,\cdots,N-1)$. 
Then one gets appropriate  factors following the rule
\be
&&\Phi_{N}(z_1)\cdots \Phi_N(z_m)=:\Phi_{N}(z_1)\cdots \Phi_N(z_m):\prod_{1\leq k<l\leq m}<\Phi_{N}(z_k)\Phi_N(z_l)>,\\
&&F_l(t_1)\cdots F_l(t_n)=:F_l(t_1)\cdots F_l(t_n):\prod_{1\leq a<b\leq n}<F_l(t_a)F_l(t_b)>.
\en
The OPE coefficients are given as follows. 
\be
<\Phi_N(z_k)\Phi_N(z_l)>&=&
\frac{\{q^2z_l/z_k\}\{pq^{2N}q^{-2}z_l/z_k\}}{\{p z_l/z_k\}\{q^{2N}z_l/z_k\}}
\nn\\
&=&(-)^{(N-1)/N}z_k^{r^*(N-1)/rN}\frac{\Gamma(q^2z_l/z_k;p,q^{2N})}{\Gamma(q^{2N}z_l/z_k;p,q^{2N})}
<\Phi_N(z_k)\Phi_N(z_l)>^{Sym},
\\
<F_l(t_a)F_l(t_b)>
&=&t_a^{2(r-1)/r}\frac{(q^2t_b/t_a;p)_\infty(t_b/t_a;p)_\infty}{(pt_b/t_a;p)_\infty(pq^{-2}t_b/t_a;p)_\infty}\\
&=&\frac{[v_a-v_b-1]}{[v_a-v_b]}<F_l(t_a)F_l(t_b)>^{Sym} 
\en
with the symmetric parts 
\be
&&<F_l(t_a)F_l(t_b)>^{Sym}=-(qt_at_b)^{1-1/r}{\frac{(t_a/t_b;p)_\infty(t_b/t_a;p)_\infty}{(pq^{-2}t_a/t_b;p)_\infty
(pq^{-2}t_b/t_a;p)_\infty}},\\
&&<\Phi_N(z_k)\Phi_N(z_l)>^{Sym}={\frac{(q^2z_k/z_l,q^2z_l/z_k;p,q^{2N})_\infty}{(q^{2N}z_k/z_l,q^{2N}z_l/z_k;p,q^{2N})_\infty}}.
\en
\item[4.] For each $l\in \{1,\cdots,N-1\}$,  symmetrize the integration variables $t^{(l)}_1,\cdots,t^{(l)}_{\la^{(l)}}$.  
We denote this procedure by ${\rm Sym}_{t^{(l)}}$.
\end{itemize}

Applying the above procedure to \eqref{nptoperator}, we obtain the following.
\be
&&\hspace{-0.8cm} \phi_{\mu_1\cdots\mu_n}(z_1,\cdots,z_n)\nn\\
&&\hspace{-0.8cm}=\oint_{\T^{N-\mu_1}} \underline{dt_1}\cdots \oint_{\T^{N-\mu_n}} \underline{dt_n}\ 
\Phi_N(z_1)F_{N-1}(t^{(N-1)}_{1})\cdots F_{\mu_1}(t^{(\mu_1)}_1)
\cdots 
\Phi_N(z_n)F_{N-1}(t^{(N-1)}_{\la^{(N-1)}})\cdots F_{\mu_n}(t^{(\mu_n)}_{\la^{(\mu_n)}})
\\
&&
\times\varphi_{\mu_1}(z_1,t^{(N-1)}_{1},\cdots,t^{(\mu_1)}_1;\{\Pi_{\mu_1,m} {q^{-2\sum_{k=2}^n<\bep_{\mu_k},h_{\mu_1,m}>}}\})
\cdots 
\varphi_{\mu_n}(z_n,t^{(N-1)}_{\la^{(N-1)}},\cdots,t^{(\mu_n)}_{\la^{(\mu_n)}};\{\Pi_{\mu_n,m}\})\\[2mm]
&&\hspace{-0.8cm}=\oint_{\T^{M}
} \underline{dt}\ 
\Phi_N(z_1)
\cdots 
\Phi_N(z_n)F_{N-1}(t^{(N-1)}_{1})\cdots F_{N-1}(t^{(N-1)}_{\la^{(N-1)}})\cdots 
F_1(t^{(1)})\cdots F_{1}(t^{(1)}_{\la^{(1)}})
\\
&&\times\prod_{l=1}^{N-2}\prod_{a=1}^{\la^{(l)}}\prod_{b=1\atop i^{(l)}_a<i^{(l+1)}_b }^{\la^{(l+1)}}\frac{[v^{(l+1)}_b-v^{(l)}_a-1/2]}{[
v^{(l+1)}_b-v^{(l)}_a+1/2]} \prod_{a=1}^{\la^{(N-1)}}\prod_{l=2}^n \frac{[u_l-v^{(N-1)}_a-1/2]}{[u_l-v^{(N-1)}_a+1/2]}\\
&&
\times
\varphi_{\mu_1}(z_1,t^{(N-1)}_{1},\cdots,t^{(\mu_1)}_1;\{\Pi_{\mu_1,m} q^{-2\sum_{k=2}^n<\bep_{\mu_k},h_{\mu_1,m}>}\})\cdots 
\varphi_{\mu_n}(z_n,t^{(N-1)}_{\la^{(N-1)}},\cdots,t^{(\mu_n)}_{\la^{(\mu_n)}};\{\Pi_{\mu_n,m}\})
\en
\be
&&\hspace{-0.8cm}=\oint_{\T^M
} \underline{dt}\ 
:\Phi_N(z_1)
\cdots 
\Phi_N(z_n)::F_{N-1}(t^{(N-1)}_{1})\cdots F_{N-1}(t^{(N-1)}_{\la^{(N-1)}}):\cdots 
:F_1(t_1^{(1)})\cdots F_{1}(t^{(1)}_{\la^{(1)}}):
\nn\\
&&\times \prod_{1\leq k<l\leq n}<\Phi_N(z_k)\Phi_N(z_l)>
\prod_{l=1}^{N-1}\prod_{1\leq a<b\leq \la^{(l)}}<F_l(t^{(l)}_a)F_l(t^{(l)}_b)>
\nn\\
&&
\times 
\prod_{l=1}^{N-2}\prod_{a=1}^{\la^{(l)}}\prod_{b=1\atop i^{(l)}_a<i^{(l+1)}_b }^{\la^{(l+1)}}\frac{[v^{(l+1)}_b-v^{(l)}_a-1/2]}{[
v^{(l+1)}_b-v^{(l)}_a+1/2]} \prod_{a=1}^{\la^{(N-1)}}\prod_{l=2}^n \frac{[u_l-v^{(N-1)}_a-1/2]}{[u_l-v^{(N-1)}_a+1/2]}
\nn\\
&&
\times
\varphi_{\mu_1}(z_1,t^{(N-1)}_{1},\cdots,t^{(\mu_1)}_1;\{\Pi_{\mu_1,m} q^{-2\sum_{k=2}^n<\bep_{\mu_k},h_{\mu_1,m}>}\})
\cdots 
\varphi_{\mu_n}(z_n,t^{(N-1)}_{\la^{(N-1)}},\cdots,t^{(\mu_n)}_{\la^{(\mu_n)}};\{\Pi_{\mu_n,m}\}).\lb{trPsiPsi}
\\[2mm]
&&\hspace{-0.8cm}=
\oint_{\T^M
} \underline{dt}\  
{:}\Phi_N(z_1)\cdots \Phi_N(z_n){:} 
 :F_{N-1}(t^{(N-1)}_{1})\cdots F_{N-1}(t^{(N-1)}_{\la^{(N-1)}}){:}\cdots 
{:}F_1(t_1^{(1)})\cdots F_{1}(t^{(1)}_{\la^{(1)}}){:}
\\
&&
\times \prod_{1\leq k<l\leq n}<\Phi_N(z_k)\Phi_N(z_l)>\prod_{l=1}^{N-1}\prod_{1\leq a<b\leq \la^{(l)}}<F_l(t^{(l)}_a)F_l(t^{(l)}_b)>^{Sym}
\\&&
\times\ 
{\rm Sym}_{t^{(1)}}\cdots {\rm Sym}_{t^{(N-1)}}
\widetilde{U}_I(t,z,\Pi),
\en
where we set 
$M=\sum_{l=1}^{N-1}\la^{(l)}=\sum_{l=1}^{N-1}(N-l)\la_l$, 
\be
&&\underline{dt}=\prod_{j=1}^na_{\mu_j,N}\prod_{l=1}^{N-1}\prod_{a=1}^{\la^{(l)}}\frac{t^{(l)}_a}{2\pi i t^{(l)}_a}, 
\en
and 
\be
&&\widetilde{U}_I(t,z,\Pi)\\
&&\hspace{-1cm}=\prod_{l=1}^{N-1}\prod_{1\leq a<b\leq \la^{(l)}}\frac{[v^{(l)}_a-v^{(l)}_b-1]}{[v^{(l)}_a-v^{(l)}_b]}\times
\prod_{l=1}^{N-2}\prod_{a=1}^{\la^{(l)}}\prod_{b=1\atop i^{(l)}_a<i^{(l+1)}_b }^{\la^{(l+1)}}\frac{[v^{(l+1)}_b-v^{(l)}_a-1/2]}{[
v^{(l+1)}_b-v^{(l)}_a+1/2]}\times \prod_{a=1}^{\la^{(N-1)}}\prod_{k=2}^n \frac{[u_k-v^{(N-1)}_a-1/2]}{[u_k-v^{(N-1)}_a+1/2]}
\\
&&\times
\varphi_{\mu_1}(z_1,t^{(N-1)}_{1},\cdots,t^{(\mu_1)}_1;\{\Pi_{\mu_1,m} q^{-2\sum_{k=2}^n<\bep_{\mu_k},h_{\mu_1,m}>}\})
\cdots \varphi_{\mu_n}(z_n,t^{(N-1)}_{\la^{(N-1)}},\cdots,t^{(\mu_n)}_{\la^{(\mu_n)}};\{\Pi_{\mu_n,m}\})
\en
with $t=(t^{(1)}_1,\cdots,t^{(1)}_{\la^{(1)}},\cdots,t^{(1)}_1,\cdots,t^{(N-1)}_{\la^{(N-1)}})$, $z=(z_1,\cdots,z_n)$, 
$z_k=q^{2u_k}\ (k=1,\cdots,n)$, $t^{(l)}_a=q^{2v^{(l)}_a}\ (l=1,\cdots,N-1, a=1,\cdots,\la^{(l)})$.  

Substituting \eqref{defvarphi} into $\widetilde{U}_I(t,z,\Pi)$, we obtain
\begin{prop}
\be
\widetilde{U}_I(t,z,\Pi)
&=&\prod_{l=1}^{N-1}\prod_{a=1}^{\la^{(l)}}\left(\frac{[v^{(l+1)}_b-v^{(l)}_a+(P+h)_{\mu_s,l+1}-C_{\mu_s,l+1}-\frac{1}{2}][1]}{[v^{(l+1)}_b-v^{(l)}_a+\frac{1}{2}][(P+h)_{\mu_s,l+1}-C_{\mu_s,l+1}(s)]}\right|_{i^{(N)}_s=i^{(l+1)}_b=i^{(l)}_a}\\[2mm]
&&\left.\qquad\qquad\qquad\times \prod_{b=1\atop i^{(l+1)}_b>i^{(l)}_a}^{\la^{(l+1)}}\frac{[v^{(l+1)}_b-v^{(l)}_a-\frac{1}{2}]}{[v^{(l+1)}_b-v^{(l)}_a+\frac{1}{2}]}\prod_{b=a+1}^{\la^{(l)}}\frac{[v^{(l)}_a-v^{(l)}_b-1]}{[v^{(l)}_a-v^{(l)}_b]}\right),
\en
where  we set $v_s^{(N)}=u_s$ and $C_{\mu_s,l+1}(s):=\sum_{j=s+1}^n<\bep_{\mu_j},h_{\mu_s,l+1}>$. 
\end{prop}

We thus obtain the following formula. 
\begin{thm}\lb{prodPhi}
For $|p|<|z_1|,\cdots, |z_n|<1$, we have
\be
&&\phi_{\mu_1\cdots\mu_n}(z_1,\cdots, z_n)=\oint_{\T^M
} \underline{dt}\  \widetilde{\Phi}(t,z)\omega_{\mu_1\cdots\mu_n}(t,z,\Pi),
\en
where  we set
\bea
&&\widetilde{\Phi}(t,z)= :\Phi_N(z_1)
\cdots 
\Phi_N(z_n):
:F_{N-1}(t^{(N-1)}_{1})\cdots F_{N-1}(t^{(N-1)}_{\la^{(N-1)}}):\cdots  :F_1(t_1^{(1)})\cdots F_{1}(t^{(1)}_{\la^{(1)}}):
\nn\\
&&\quad
\times \prod_{1\leq k<l\leq n}<\Phi_N(z_k)\Phi_N(z_l)>^{Sym}
\prod_{l=1}^{N-1}\prod_{1\leq a<b\leq \la^{(l)}}<F_l(t^{(l)}_a)F_l(t^{(l)}_b)>^{Sym},
\lb{Phitilde}\\
&&
\omega_{\mu_1\cdots\mu_n}(t,z,\Pi)= \mu^+(z)\ 
\widetilde{W}_I(t,z,\Pi),\lb{def:omega}\\
&&\mu^+(z)=\prod_{1\leq k<l\leq n}(-)^{(N-1)/N}z_k^{r^*(N-1)/rN}\frac{\Gamma(q^2z_k/z_l;p,q^{2N})}{\Gamma(q^{2N}z_k/z_l;p,q^{2N})},
\lb{def:mup}\\
&&
\widetilde{W}_I(t,z,\Pi)= {\rm Sym}_{t^{(1)}}\cdots {\rm Sym}_{t^{(N-1)}}
\widetilde{U}_I(t,z,\Pi).\lb{def:Wtilde}
\ena
\end{thm}
Note that  $\widetilde{\Phi}(t,z)$ is a symmetric function in $z_1,\cdots,z_n$ as well as in $t^{(l)}_1,\cdots,t^{(l)}_{\la^{(l)}}$ for each $l\in \{1,\cdots,N-1\}$. 
The function $\widetilde{W}_I(t,z,\Pi)$ is the  weight function of type $\sln$, 
which is an elliptic and dynamical analogue of the trigonometric 
one in \cite{Mimachi,MN,RTV}. 

\noindent
{\it Remark.}\  In our derivation, the weight functions are determined modulo adding co-boundary terms associated with the multiple integral. 
In \eqref{def:Wtilde} we have resolved this ambiguity by requiring that their trigonometric limit coincide with the known results  in \cite{Mimachi,MN,RTV}. An alternative and more intrinsic way to fix this ambiguity is provably to require the triangular property in Proposition \ref{triangular}. 
We apply the same argument to the proof of Propositions \ref{transitionProp} and \ref{shuffleprod}.

For a partition $I=(I_1,\cdots,I_N)\in \cI_\la$, 
let  $I_k=\{i_{k,1}<\cdots<i_{k,\la_k}\}$ $(k=1,\cdots,N)$.
The dynamical shift term $C_{\mu_s,l+1}$ appearing in $\widetilde{U}_I(t,z,\Pi)$ has the following combinatorial expression. 
\begin{prop}\lb{combC}
\be
&&
C_{\mu_s,l+1}(s)
=\left\{\mmatrix{\la_{\mu_s}-\la_{l+1}-\tilde{s}+m_{\mu_s,l+1}(s)-1&\quad \mbox{if }\ s<i_{l+1,\la_{l+1}}\cr
\la_{\mu_s}-\tilde{s}\qquad\qquad &\quad \mbox{if }\ s>i_{l+1,\la_{l+1}}\cr}\right.
\en
where for $s\in [1,n]$ we define $\tilde{s}$ by $i_{\mu_s,\tilde{s}}=s$  and $m_{\mu_s,l+1}(s)$  by  
\be
&&m_{\mu_s,l+1}(s)=\mbox{min}\{1\leq j\leq \la_{l+1}\ | \ s<i_{l+1,j}\ \} \quad \mbox{for }\ s<i_{l+1,\la_{l+1}}.
\en
\end{prop}
\noindent
{\it Proof.}\ 
Note by definition $\mu_s\leq l$ and 
\be
&&\sum_{j=s+1}^n<\bep_{\mu_j},h_{\mu_s,l+1}>=\sum_{j=s+1}^{n}(\delta_{\mu_j,\mu_s}-\delta_{\mu_j,l+1}). 
\en
Let  $I_{\mu_s}=\{i_{\mu_s,1}<\cdots<s=i_{\mu_s,\tilde{s}}<\cdots<i_{\mu_s,\la_{\mu_s}}\}$, then 
$\sum_{j=s+1}^{n}\delta_{\mu_j,\mu_s}=\la_{\mu_s}-\tilde{s}$. If $s<i_{l+1,\la_{l+1}}$  there exists 
$m_{\mu_s,l+1}(s)$ so that $s<i_{l+1,j}\in I_{l+1}\ (m_{\mu_s,l+1}(s)\leq j\leq\la_{l+1})$. Hence 
$\sum_{j=s+1}^{n}\delta_{\mu_j,l+1}=\la_{l+1}-m_{\mu_s,l+1}(s)+1$. If $s>i_{l+1.\la_{l+1}}$ then 
$\sum_{j=s+1}^{n}\delta_{\mu_j,l+1}=0$. 
\qed


\subsection{Entire function version}
Let us set 
\bea
&&H_\la(t,z):=\prod_{l=1}^{N-1}\prod_{a=1}^{\la^{(l)}}\prod_{b=1}^{\la^{(l+1)}}\left[v^{(l+1)}_b-v^{(l)}_a+\frac{1}{2}\right].\lb{Efunc}
\ena
The following gives an entire function version of $\widetilde{W}_I$
\bea
&&W_I(t,z,\Pi)=H_\la(t,z)\widetilde{W}_I(t,z,\Pi)
= {\rm Sym}_{t^{(1)}}\cdots {\rm Sym}_{t^{(N-1)}}
{U}_I(t,z,\Pi),\lb{entireW}
\ena
where
\bea
&&{U}_I(t,z,\Pi)= 
\prod_{l=1}^{N-1}\prod_{a=1}^{\la^{(l)}}\left(\left.
\frac{\left[v^{(l+1)}_b-v^{(l)}_a+(P+h)_{\mu_s,l+1}-C_{\mu_s,l+1}(s)-\frac{1}{2}\right][1]}{[(P+h)_{\mu_s,l+1}-C_{\mu_s,l+1}(s)]}
\right|_{i^{(N)}_s=i^{(l+1)}_b=i^{(l)}_a}\right.\nn\\
&&\left.\qquad\qquad\times \prod_{b=1\atop i^{(l+1)}_b>i^{(l)}_a}^{\la^{(l+1)}}{\left[v^{(l+1)}_b-v^{(l)}_a-\frac{1}{2}\right]}
\prod_{b=1\atop i^{(l+1)}_b<i^{(l)}_a}^{\la^{(l+1)}}{\left[v^{(l+1)}_b-v^{(l)}_a+\frac{1}{2}\right]}
\prod_{b=a+1}^{\la^{(l)}}\frac{[v^{(l)}_b-v^{(l)}_a+1]}{[v^{(l)}_b-v^{(l)}_a]}\right).\nn\\
&&\lb{entireU}
\ena
In the trigonometric $(p\to 0)$ and non-dynamical (neglecting the factors depending on $P+h$ ) limit the $W_I$ coincides  with those 
discussed in \cite{RTV} by making the shift $v^{(l)} \mapsto v^{(l)}+l/2$. See Sec.\ref{modW}. Note however that in this limit our $\tW_I$ are slightly different from those 
in (6.2) \cite{RTV} due to the difference between our $H_\la(t,z)$ and $E(t,h)$ in   (6.1) from \cite{RTV}.

\subsection{Combinatorial description } 
One can determine each factor in ${U}_I(t,z,\Pi)$ by using the following  one to one correspondence to a pair of variables $(t^{(l)}_a,t^{(l')}_b)$ in a table introduced in Sec.\ref{combNot} \cite{GRTV,RTV}.  As above  we set $t^{(N)}_a=z_a$. 
 \begin{itemize}
 \item for each pair $(t^{(l)}_a,t^{(l+1)}_b)$ $(l=1,\cdots,N-1)$ in the same row (i.e. the $i^{(l)}_a=i^{(l+1)}_b$-th row), assign a factor 
 \be
 &&\frac{\left[v^{(l+1)}_b-v^{(l)}_a+(P+h)_{\mu_s,l+1}-C_{\mu_s,l+1}(s)-\frac{1}{2}\right][1]}{[(P+h)_{\mu_s,l+1}-C_{\mu_s,l+1}(s)]}, 
 \en
 where $s$ is the index of $z_s$ located in the same row ( $i^{(l)}_a=i^{(l+1)}_b=i^{(N)}_s$). 
\item  for each pair $(t^{(l)}_a,t^{(l+1)}_b)$ $(l=1,\cdots,N-1)$, where  $t^{(l+1)}_b$ is located in the lower row than $  t^{(l)}_a$ ( i.e. $i^{(l)}_a<i^{(l+1)}_b$), 
assign 
\be
&&\left[v^{(l+1)}_b-v^{(l)}_a-\frac{1}{2}\right].
\en
\item  for each pair $(t^{(l)}_a,t^{(l+1)}_b)$ $(l=1,\cdots,N-1)$, where  $t^{(l+1)}_b$ is located in the upper row than $  t^{(l)}_a$ (i.e. $i^{(l)}_a>i^{(l+1)}_b$), 
assign 
\be
&&\left[v^{(l+1)}_b-v^{(l)}_a+\frac{1}{2}\right].
\en
\item for each pair $(t^{(l)}_a,t^{(l)}_b)$ $(1\leq a<b\leq \la^{(l)})$ in the $l$-th column $(l=1,\cdots,N-1)$, assign 
\be
&&\frac{[v^{(l)}_b-v^{(l)}_a+1]}{[v^{(l)}_b-v^{(l)}_a]}.
\en
 \end{itemize}
\noindent
{\it Example 3.}  For the case $N=3, n=4, I=\{I_1=\{2\},I_2=\{1,4\},I_3=\{3\}\}=I_{2132}$\\[2mm]
$$
\begin{array}{|c|c|c|}\hline
              &t^{(2)}_1&z_1\\ \hline
t^{(1)}_1&t^{(2)}_2&z_2\\ \hline
              &             &z_3\\ \hline
            &t^{(2)}_3&z_4\\ \hline
 \end{array}       
$$
we obtain
\be
&&U_I(t,z,\Pi)\\
&&=\frac{[v^{(2)}_2-v^{(1)}_1+(P+h)_{1,2}-C_{1,2}(2)-1/2][1]}{[(P+h)_{1,2}-C_{1,2}(2)]}
\frac{[u_1-v^{(2)}_1+(P+h)_{2,3}-C_{2,3}(1)-1/2][1]}{[(P+h)_{2,3}-C_{2,3}(1)]}\\
&&\times\frac{[u_2-v^{(2)}_2+(P+h)_{1,3}-C_{1,3}(2)-1/2][1]}{[(P+h)_{1,3}-C_{1,3}(2)]}
\frac{[u_4-v^{(2)}_3+(P+h)_{2,3}-C_{2,3}(4)-1/2][1]}{[(P+h)_{2,3}-C_{2,3}(4)]}\\
&&\hspace{-0cm}\times[v^{(2)}_3-v^{(1)}_1-1/2][u_2-v^{(2)}_1-1/2][u_3-v^{(2)}_1-1/2][u_4-v^{(2)}_1-1/2][u_3-v^{(2)}_2-1/2][u_4-v^{(2)}_2-1/2]\\
&&\times [v^{(2)}_1-v^{(1)}_1+1/2][u_1-v^{(2)}_2+1/2][u_1-v^{(2)}_3+1/2][u_2-v^{(2)}_3+1/2][u_3-v^{(2)}_3+1/2]\\
&&\prod_{1\leq a<b\leq 3}\frac{[v^{(2)}_b-v^{(2)}_a+1]}{[v^{(2)}_b-v^{(2)}_a]}.
\en

Furthermore, 
in order to compute $C_{\mu_s,l+1}(s)$ by Proposition \ref{combC}, it is convenient to draw a $n\times N$ table of the elements in $I_l=\{i_{l,1}< \cdots < i_{l,\la_l} \}$\ $(l=1,\cdots,N)$ 
obtained by the rule: put $i_{l,a}$ into the $(i_{l,a}, l)$ box. 

\noindent
{\it Example 4.}\ For the case in {\it Example 3}, we have
$$
\begin{array}{|c|c|c|}\hline
              &i_{2,1}&      \\ \hline
i_{1,1}&             &       \\ \hline
              &             &i_{3,1}\\ \hline
            & i_{2,2}           &  \\ \hline
\end{array}       
$$     
For $C_{1,2}(2)$, we have $s=i_{1,1}=2<i_{2,2}=4$, $m_{1,2}(2)=2$. Hence 
\be
C_{1,2}(2)=1-2-1+2-1=-1=\sum_{j=3,4}<\bep_{\mu_j},h_{1,2}>.
\en
For $C_{2,3}(1)$, we have $s=i_{2,1}=1<i_{3,1}=3$, $m_{2,3}(1)=1$. Hence 
\be
C_{2,3}(1)=2-1-1+1-1=0=\sum_{j=2,3,4}<\bep_{\mu_j},h_{2,3}>.
\en
For $C_{1,3}(2)$, we have $s=i_{1,1}=2<i_{3,1}=3$, $m_{1,3}(2)=1$. Hence 
\be
C_{1,3}(2)=1-1-1+1-1=-1=\sum_{j=3,4}<\bep_{\mu_j},h_{1,3}>.
\en
For $C_{2,3}(4)$, we have $s=i_{2,2}=4>i_{3,1}=3$. Hence 
\be
C_{2,3}(4)=2-2=0. 
\en

\section{Properties of the Elliptic Weight Functions}\lb{modW}
In this section we discuss some basic properties of  the elliptic weight functions. 
We consider the following weight functions obtained from those in the last section by 
$v^{(l)}_a \mapsto v^{(l)}_a+l/2$ $(l=1,\cdots,N, a=1,\cdots, \la^{(l)})$.
\bea
&&W_I(t,z,\Pi)=  {\rm Sym}_{t^{(1)}}\cdots {\rm Sym}_{t^{(N-1)}}
{U}_I(t,z,\Pi),\lb{entireWmod}\\[2mm]
&&{U}_I(t,z,\Pi)=\prod_{l=1}^{N-1}\prod_{a=1}^{\la^{(l)}}\left(\left.
\frac{\left[v^{(l+1)}_b-v^{(l)}_a+(P+h)_{\mu_s,l+1}-C_{\mu_s,l+1}(s)\right][1]}{[(P+h)_{\mu_s,l+1}-C_{\mu_s,l+1}(s)]}
\right|_{i^{(N)}_s=i^{(l+1)}_b=i^{(l)}_a}\right.\nn\\
&&\left.\qquad\qquad\times \prod_{b=1\atop i^{(l+1)}_b>i^{(l)}_a}^{\la^{(l+1)}}{\left[v^{(l+1)}_b-v^{(l)}_a\right]}
\prod_{b=1\atop i^{(l+1)}_b<i^{(l)}_a}^{\la^{(l+1)}}{\left[v^{(l+1)}_b-v^{(l)}_a+1\right]}
\prod_{b=a+1}^{\la^{(l)}}\frac{[v^{(l)}_b-v^{(l)}_a+1]}{[v^{(l)}_b-v^{(l)}_a]}\right).\nn\\
&&\lb{entireUmod}
\ena
Accordingly using  
\bea
&&H_\la(t,z):=\prod_{l=1}^{N-1}\prod_{a=1}^{\la^{(l)}}\prod_{b=1}^{\la^{(l+1)}}\left[v^{(l+1)}_b-v^{(l)}_a+1\right],
\lb{Hfuncmod}
\ena
the meromorphic  function version of ${W}_I$ is given by 
\bea
&&\widetilde{W}_I(t,z,\Pi)=\frac{W_I(t,z,\Pi)}{H_\la(t,z)}. \lb{meromWmod}
\ena
   
\subsection{Triangular property}\lb{SecTriangular}
For $I,J\in \cI_{\la}$, let $I^{(l)}=\{i^{(l)}_1<\cdots<i^{(l)}_{\la^{(l)}}\}$ and $J^{(l)}=\{j^{(l)}_1<\cdots<j^{(l)}_{\la^{(l)}}\}$ 
$(l=1,\cdots,N)$.
 Define a partial ordering $\leqslant$ by 
\be
I\leqslant J \Leftrightarrow i^{(l)}_a \leq j^{(l)}_a\qquad \forall l, a. 
\en
Let us denote by $t=z_I$ the specialization $t^{(l)}_a=z_{i^{(l)}_a}$ $(l=1,\cdots,N-1, a=1,\cdots,\la^{(l)})$\cite{RTV}.
The weight function has the following triangular property.
\begin{prop}\lb{triangular}
For $I,J\in \cI_\la$, 
\begin{itemize}
\item[(1)] $\tW_{J}(z_I,z,\Pi)=0$ unless $I\leqslant J$.
\item[(2)] 
\be
&&\tW_{I}(z_I,z,\Pi)=\prod_{1\leq k<l\leq N}\prod_{a\in I_k}\prod_{b\in I_l\atop a<b}\frac{[u_b-u_a]}{[u_b-u_a+1]}.
\en
\end{itemize}

\end{prop}
\noindent
{\it Proof.}\ (1) follows from the same argument as in \cite{RTV} Lemma 6.2. \\
(2) follows from  
\be
W_I(z_I,z,\Pi)&=&\prod_{1\leq k<l\leq N}\left[\prod_{a\in I_k}\left(\prod_{b\in  I_k}[u_b-u_a+1]
\prod_{b\in I_l\atop a<b}[u_b-u_a]\prod_{b\in I_l\atop a>b}[u_b-u_a+1]\right)\right]\\
&&\times
\prod_{1\leq k<l\leq N-1}\prod_{a\in I_k}\prod_{b\in I_l}[u_b-u_a+1][u_a-u_b+1],\\
H_\la(z_I,z)&=&\prod_{1\leq k<l\leq N}\left[\prod_{a\in I_k}\left(\prod_{b\in  I_k}[u_b-u_a+1]
\prod_{b\in I_l}[u_b-u_a+1]\right)\right]\\
&&\times
\prod_{1\leq k<l\leq N-1}\prod_{a\in I_k}\prod_{b\in I_l}[u_b-u_a+1][u_a-u_b+1].
\en
\qed

For $\sigma\in \gS_n$, let us denote $\sigma^{-1}(I)=I_{\mu_{\sigma(1)}\cdots\mu_{\sigma(n)}}$ and 
$\sigma(z)=(z_{\sigma(1)},\cdots,z_{\sigma(n)})$.  
Following  \cite{RTV}, 
let us set $\widetilde{W}_{\sigma,I}(t,z,\Pi)=\widetilde{W}_{\sigma^{-1}(I)}(t,\sigma(z),\Pi)$ 
and $\tW_{\id,I}(t,z\Pi)=\tW_I(t,z,\Pi)$. 
Let us  consider the matrix $\hW_{\sigma}(z,\Pi)$, whose $(I,J)$ element is given by $\tW_{\sigma,J}(z_I,z,\Pi)$ $(I,J\in \cI_\la)$. 
We put the matrix elements in the decreasing order with respect to $\leqslant$. 
Then Proposition \ref{triangular} yields that the matrix $\hW_{\id}(z,\Pi)$ is lower triangular, whereas 
$\hW_{\sigma_0}(z,\Pi)$ with $\sigma_0$ being the longest element in $\gS_n$ is upper triangular. 
In particular, for generic $u_a\ (a=1,\cdots,n)$,  $\hW_{\sigma}(z,\Pi)$ is invertible.

\subsection{Transition property}\lb{TransProp}
\begin{prop}\lb{transitionProp}
Let $I=I_{\mu_1\cdots\mu_i\mu_{i+1}\cdots\mu_n}\in \cI_\la$. 
\bea
&&
\widetilde{W}
_{I_{\cdots\ \mu_{i+1}\mu_{i}\cdots}}
(t, 
\cdots,z_{i+1},z_{i},\cdots
,\Pi)\nn\\
&&=
\sum_{\mu_{i}',\mu_{i+1}'}\bR(z_{i}/z_{i+1},\Pi q^{-2\sum_{j=i}^{n}<\bep_{\mu_j},h>})_{\mu_{i}\mu_{i+1}}^{\mu_{i}'\mu_{i+1}'}\ 
\widetilde{W}_{I_{\cdots\ \mu'_i\mu'_{i+1}\cdots}}(t, 
\cdots,z_i,z_{i+1},\cdots
,\Pi).\lb{transsi}
\ena
\end{prop}
\noindent
{\it Proof.}\ 
From Theorem \ref{prodPhi} we have 
\be
&& \phi_{\cdots \mu_{i+1} \mu_i\cdots}(\cdots, z_{i+1}, z_i, \cdots)=
\oint_{\T^M} \underline{d t}\ {\widetilde{\Phi}(t,z)}\ {\omega}_{\cdots\mu_{i+1} \mu_i\cdots}(t,\cdots, z_{i+1}, z_i, \cdots;\Pi)
\en
where we used the symmetry property of ${\widetilde{\Phi}(t,z)}$ under any permutations of $z_1,\cdots, z_n$.  
Using the exchange relation \eqref{typeIcr}  in the left hand side, we obtain
\be
&&\phi_{\cdots \mu_{i+1} \mu_i\cdots}(\cdots, z_{i+1}, z_i, \cdots)\\
&&=\sum_{\mu'_i,\mu'_{i+1}}R(z_{i}/z_{i+1},\Pi q^{2\sum_{j=1}^{i-1}<\bep_{\mu_j},h>})_{\mu_i\mu_{i+1}}^{\mu'_i\mu'_{i+1}}\ \phi_{\cdots \mu'_{i} \mu'_{i+1}\cdots}(\cdots, z_i, z_{i+1}, \cdots)\\
&&=\oint_{\T^M} \underline{d t}\ {\widetilde{\Phi}(t,z)}\ 
\sum_{\mu'_i,\mu'_{i+1}}R(z_{i}/z_{i+1},\Pi q^{-2\sum_{j=i}^{n}<\bep_{\mu_j},h>})_{\mu_i\mu_{i+1}}^{\mu'_i\mu'_{i+1}}\
{\omega}_{\cdots \mu'_i\mu'_{i+1}\cdots}(t,\cdots,  z_i,z_{i+1}, \cdots;\Pi).
\en
Note that $\mu(z)$ \eqref{def:mu} in the $R$ matrix is related to $\mu^+(z)$ \eqref{def:mup} by
\be
&&\mu(z_i/z_{i+1})=\mu^+(\cdots,z_{i+1},z_i,\cdots )/\mu^+(\cdots,z_{i},z_{i+1},\cdots ).
\en
Comparing the integrand we obtain the desired relation.
\qed

Let $s_i=(i\ i+1)\in \gS_n (i=1,\cdots,n-1)$ denote the adjacent transpositions. 
Let us set
\be
&&\cR^{(s_i,\id)}(z,\Pi)_{I}^{I}=\bR(z_{i}/z_{i+1},\Pi q^{2\sum_{j=1}^{i-1}<\bep_{\mu_j},h>})_{\mu_{i}\mu_{i+1}}^{\mu_{i}\mu_{i+1}},\\
&&\cR^{(s_i,\id)}(z,\Pi)_I^{s_i(I)}=\bR(z_{i}/z_{i+1},\Pi q^{2\sum_{j=i}^{i-1}<\bep_{\mu_j},h>})_{\mu_{i}\mu_{i+1}}^{\mu_{i+1}\mu_{i}}.
\en
Then one can rewrite \eqref{transsi} as
\be
&&
\tW_{s_i,I}(t,z,\Pi)\nn\\
&&=\left\{\mmatrix{
\tW_{I}(t,z,\Pi)\qquad\qquad\qquad\qquad\qquad&\mbox{if}\ s_i(I)=I\cr
\cR^{(s_i,\id)}(z,\Pi q^{-2\sum_{j=1}^{n}<\bep_{\mu_j},h>})_{I}^{I}
\tW_{\id,I}(t,z,\Pi)\qquad\qquad\\
\qquad\qquad\qquad+\cR^{(s_i,\id)}(z,\Pi q^{-2\sum_{j=1}^{n}<\bep_{\mu_j},h>})_I^{s_i(I)}
\tW_{\id,s_i(I)}(t,z,\Pi)&       \mbox{if}\ s_i(I)\not=I\cr  }\right..
\en

In general, let us consider the $n$-point operator
\be
&&\phi_{\sigma^{-1}(I)}(\sigma(z))=\Phi_{\mu_{\sigma(1)}}(z_{\sigma(1)})\cdots \Phi_{\mu_{\sigma(n)}}(z_{\sigma(n)}). 
\en
By using the exchange relation \eqref{typeIcr} repeatedly we obtain
\be
&&\phi_{\sigma^{-1}(I)}(\sigma(z))=\sum_{I'}\cR^{(\sigma,\sigma')}(z,\Pi)_{I}^{I'}\phi_{{\sigma'}^{-1}(I')}(\sigma'(z))
\en
where $\cR^{(\sigma,\sigma')}(z,\Pi)_{I}^{I'}$ denotes a coefficient given by a certain sum of products of  the $R$ matrix elements in \eqref{Relements}. 
Then by the same argument as in the proof of Proposition \ref{transitionProp} we obtain 
\bea
&&
\tW_{\sigma,I}(t,z,\Pi)=\sum_{I'}\cR^{(\sigma,\sigma')}(z,\Pi 
q^{-2\sum_{j=1}^{n}<\bep_{\mu_j},h>}
)_{I}^{I'}
\tW_{\sigma', I'}(t,z,\Pi).\lb{WRW}
\ena
Let us define the matrix $\cR^{(\sigma,\sigma')}(z,\Pi)$ by
\be
&&\cR^{(\sigma,\sigma')}(z,\Pi)=\left(\ \cR^{(\sigma,\sigma')}(z,\Pi)_{I}^{J} \ \right)_{I,J\in \cI_\la}.
\en
Then we can rewrite \eqref{WRW}
\bea
&&
\hW_{\sigma}(z,\Pi)=
\hW_{\sigma'}(z,\Pi)\ {}^t\cR^{(\sigma,\sigma')}(z,\Pi q^{-2\sum_{j=1}^{n}<\bep_{\mu_j},h>}).\lb{WWtR}
\ena
or
\bea
 &&\hW_{\sigma'}(z,\Pi)^{-1}\ \hW_{\sigma}(z,\Pi)= {}^t\cR^{(\sigma,\sigma')}(z,\Pi q^{-2\sum_{j=1}^{n}<\bep_{\mu_j},h>}).\lb{WinvWtR}
\ena
By taking the transposition and the shift $\Pi\mapsto \Pi q^{2\sum_{j=1}^{n}<\bep_{\mu_j},h>}$, we obtain
\bea
&&{}^t\hW_{\sigma}(z,\Pi q^{2\sum_{j=1}^{n}<\bep_{\mu_j},h>})\left({}^t\hW_{\sigma'}(z,\Pi q^{2\sum_{j=1}^{n}<\bep_{\mu_j},h>})\right)^{-1}= \cR^{(\sigma,\sigma')}(z,\Pi)\ .\lb{tWRtW}
\ena
In addition,  we have 
\begin{prop}
\bea
&&{}^t\cR^{(\sigma,\sigma')}(z,\Pi q^{-2\sum_{j=1}^{n}<\bep_{\mu_j},h>})=\cR^{(\sigma,\sigma')}(z,\Pi^{-1}). \lb{RtR}
\ena
\end{prop}
\noindent
{\it Proof.}\ It is enough to show the case that $\cR^{(\sigma,\sigma')}(z,\Pi)$ is given by 
\be
&&R^{(23)}(z_2/z_3,\Pi)R^{(13)}(z_1/z_3,\Pi q^{2h^{(2)}})R^{(12)}(z_1/z_2,\Pi). 
\en
The desired equality follows from the properties of the $R$ matrix \eqref{def:Rmat} such as  ${}^tR(z,\Pi)=R(z,\Pi^{-1})$ and $R(z,\Pi q^{2(h^{(1)}+h^{(2)})})=R(z,\Pi)$
as well as 
the dynamical Yang-Baxter equation \eqref{DYBE}. \qed

\subsection{Orthogonality}
This is an elliptic and dynamical analogue of the same property given in \cite{RTV}.

\begin{prop}
For $J,K\in \cI_\la$, 
\be
&&\sum_{I\in \cI_\la}\frac{W_J(z_I,z,\Pi^{-1}q^{2\sum_{j=1}^{n}<\bep_{\mu_j},h>})W_{\sigma_0(K)}(z_I,\sigma_0(z),\Pi)}{Q(z_i)R(z_I)S(z_I)^2}=\delta_{J,K},
\en
where $\sum_{j=1}^n\bep_{\mu_j}$ is the weight associated with $\la$ (Sec.\ref{combNot}), and
\be
&&Q(z_I)=\prod_{1\leq k<l\leq N}\prod_{a\in I_k}\prod_{b\in I_l}[u_b-u_a+1],\\
&&R(z_I)=\prod_{1\leq k<l\leq N}\prod_{a\in I_k}\prod_{b\in I_l}[u_b-u_a],\\
&&S(z_I)=\prod_{1\leq k<l\leq N}^N\prod_{a,b\in I_k }[u_a-u_b+1]\prod_{1\leq k\not=l\leq N-1}^N\prod_{a\in I_k }\prod_{b\in I_l }[u_a-u_b+1].
\en
\end{prop}
\noindent
{\it Proof.}\ The proof is parallel to the one in \cite{RTV} except for the dynamical shift.     
From \eqref{WinvWtR}, \eqref{tWRtW} and \eqref{RtR} with $\sigma=\sigma_0$ and $\sigma'=\id$, we have
\be
&&\hW_{\id}(z,\Pi)^{-1}\hW_{\sigma_0}(z,\Pi)={}^t  \hW_{\sigma_0}(z,\Pi^{-1}q^{2\sum_{j=1}^{n}<\bep_{\mu_j},h>}) \left({}^t \hW_{\id}(z,\Pi^{-1}q^{2\sum_{j=1}^{n}<\bep_{\mu_j},h>})\right)^{-1}.
\en
Hence 
\be
&&\hW_{\id}(z,\Pi)\ {}^t  \hW_{\sigma_0}(z,\Pi^{-1}q^{2\sum_{j=1}^{n}<\bep_{\mu_j},h>})=\hW_{\sigma_0}(z,\Pi)\  {}^t \hW_{\id}(z,\Pi^{-1}q^{2\sum_{j=1}^{n}<\bep_{\mu_j},h>}).
\en
Since the LHS is a lower triangular matrix and the RHS is an upper triangular matrix, this must be a diagonal matrix. 
Let us denote it by $S$. It's diagonal entries are obtained from Proposition \ref{triangular} (2) as 
\be
&&S_{II}=\tW_{I}(z_I,z,\Pi)\tW_{\sigma_0(I)}(z_{I},\sigma_0(z),\Pi^{-1}q^{2\sum_{j=1}^{n}<\bep_{\mu_j},h>})=\frac{R(z_I)}{Q(z_I)}.
\en    
We then obtain 
\be
&&{}^t\hW_{\id}(z,\Pi^{-1}q^{2\sum_{j=1}^{n}<\bep_{\mu_j},h>})S^{-1}\hW_{\sigma_0}(z,\Pi)=\id
\en
Taking the $(J,K)$ component of this relation, we obtain
\be
&&\sum_{I\in \cI_\la}\tW_J(z_I,z,\Pi^{-1}q^{2\sum_{j=1}^{n}<\bep_{\mu_j},h>})\frac{Q(z_I)}{R(z_I)}\tW_{\sigma_0(K)}(z_I,\sigma_0(z),\Pi)=\delta_{J,K}
\en 
Furthermore noting 
\be
&&\tW_{J}(z_I,z,\Pi)=\frac{W_{J}(z_I,z,\Pi)}{H_\la(z_I,z)} \qquad\qquad\qquad \forall J\in \cI_\la
\en
and from the proof of Proposition \ref{triangular}
\be
&&H_\la(z_I,z)=H_\la(z_I,\sigma_0(z))=Q(z_I) S(z_I)
\en
we obtain the desired result.
\qed

\noindent
{\it Remark.} \ In \cite{RTV17}, the following function is used instead of our $H_\la(t,z)$. 
\be
&&E_\la(t,z)=\prod_{l=1}^{N-1}\prod_{a=1}^{\la^{(l)}}\prod_{b=1}^{\la^{(l)}}[v^{(l)}_b-v^{(l)}_a+1].
\en
The specialization $E_\la(z_I,z)$ coincides with our $S(z_I)$. 

\subsection{Quasi-periodicity}
Remember  that we set $t^{(l)}_a=q^{2v^{(l)}_a}$, $z_k=q^{2u_k}$ and $\Pi_{j,k}=q^{2(P+h)_{j,k}}$.  Note that  $t^{(l)}_a\mapsto pt^{(l)}_a\Leftrightarrow v^{(l)}_a\mapsto v^{(l)}_a+r$ and
$t^{(l)}_a\mapsto e^{-2\pi i}t^{(l)}_a\Leftrightarrow v^{(l)}_a\mapsto v^{(l)}_a+r\tau$. 
From \eqref{thetaquasiperiod} and Proposition \ref{combC} we obtain the following statement.   
\begin{prop}\lb{quasiperiod}
For $I\in \cI_\la$, the weight functions $W_I(t,z,\Pi)$ have the following quasi-periodicity. 
\be
&&W_I(\cdots,pt^{(l)}_a,\cdots, z,\Pi)=(-1)^{\la^{(l+1)}+\la^{(l-1)}}W_I(\cdots,t^{(l)}_a,\cdots, z,\Pi),\\
&&W_I(\cdots,e^{-2\pi i}t^{(l)}_a,\cdots, z,\Pi)\\
&&=(-e^{-\pi i \tau})^{\la^{(l+1)}+\la^{(l-1)}}\\
&&\qquad\times \exp\left\{-\frac{2\pi i}{r}\left((\la^{(l+1)}+\la^{(l-1)})v^{(l)}_a-\sum_{b=1}^{\la^{(l+1)}}v^{(l+1)}_b-\sum_{c=1}^{\la^{(l-1)}}v^{(l-1)}_c-(P+h)_{l,l+1}-\la_{l+1}\right)\right\}\\
&&\qquad \times W_I(\cdots,t^{(l)}_a,\cdots, z,\Pi) \qquad\qquad (1\leq a\leq \la^{(l)}, 1\leq l\leq N-1).
\en
\end{prop}

Furthermore noting
\be
&&H_\la(\cdots,pt^{(l)}_a,\cdots ,z)=(-1)^{\la^{(l+1)}+\la^{(l-1)}}H_\la(\cdots,t^{(l)}a,\cdots ,z),\\
&&H_\la(\cdots,e^{-2\pi i}t^{(l)}_a,\cdots ,z)=(-e^{-\pi i \tau})^{\la^{(l+1)}+\la^{(l-1)}}\\
&&\qquad\times \exp\left\{-\frac{2\pi i}{r}\left((\la^{(l+1)}+\la^{(l-1)})v^{(l)}_a-\sum_{b=1}^{\la^{(l+1)}}v^{(l+1)}_b-\sum_{c=1}^{\la^{(l-1)}}v^{(l-1)}_c
-\la_{l+1}-\la_{l}\right)\right\}\\
&&\qquad \times H_\la(\cdots,t^{(l)}_a,\cdots, z),
\en
we obtain
\begin{prop}\lb{meroW}
The meromorphic weight functions $\widetilde{W}_I(t,z,\Pi)$ have the following quasi-periodicity.
\be
&&\widetilde{W}_I(\cdots,pt^{(l)}_a,\cdots, z,\Pi)=\widetilde{W}_I(\cdots,t^{(l)}_a,\cdots, z,\Pi),\\
&&\widetilde{W}_I(\cdots,e^{-2\pi i}t^{(l)}_a,\cdots, z,\Pi)=\exp\left\{\frac{2\pi i}{r}\left((P+h)_{l,l+1}-\la_{l}\right)\right\}
\widetilde{W}_I(\cdots,t^{(l)}_a,\cdots, z,\Pi) \nn\\
&&\qquad\qquad (l=1,\cdots,N-1,\ a=1, \cdots,\la^{(l)}).
\en
\end{prop}

For $\la=(\la_1,\cdots,\la_N)\in \N^N$, $|\la|=n$, let $z^{(n)}=(z_1,\cdots,z_n)\in (\C^*)^n$.   
For $I=I_{\mu_1\cdots \mu_n}\in \cI_{\la}$, we denote by $\Pi_I$ a set of dynamical parameters $\{\Pi_{\mu_k,j}=q^{2(P+h)_{\mu_k,j}}\ (k=1,\cdots, n, j=\mu_k+1,\cdots,N) \}$,where  ${(P+h)_{j,k}}\in \C/r\Z\ (1\leq j<k\leq N)$, and set $\Pi_\la=\cup_{I\in\cI_\la}\Pi_I$. 
\begin{dfn}
For $\la=(\la_1,\cdots,\la_N)\in \N^N$, $|\la|=n$, 
we define $\cM^{(n)}_\la(z^{(n)},\Pi_\la)$ to be  the space of meromorphic functions $F(t;z,\Pi)$ of $M=\sum_{l=1}^{N-1}(N-l)\la_l$ variables  $t=(t^{(1)}_1,\cdots,t^{(1)}_{\la^{(1)}},\cdots,$ \\ $ t^{(N-1)}_{1},\cdots,t^{(N-1)}_{\la^{(N-1)}})$ such that 
\begin{itemize}
\item[(1)] $F(t;z,\Pi)$ is symmetric in $t^{(l)}_{1},\cdots,t^{(l)}_{\la^{(l)}}$ for each $l\in \{1,\cdots,N-1\}$. 
\item[(2)] ${F}(t;z,\Pi)$ has the quasi-periodicity
\be
&&{F}(\cdots, pt^{(l)}_a,\cdots;z,\Pi )={F}(t;z,\Pi),\\
&&{F}(\cdots, e^{-2\pi i}t^{(l)}_a;\cdots,z,\Pi )=\exp\left\{\frac{2\pi i}{r}\left((P+h)_{l,l+1}-\la_{l}\right)\right\}{F}(t;z,\Pi)
\en
$(l=1,\cdots,N-1,\ a=1,\cdots,\la^{(l)})$.

\end{itemize}
\end{dfn}

Let us consider the subspace $\cM^{+(n)}_\la(z^{(n)},\Pi_\la):={\rm Span}_\C\{\ \widetilde{W}_I(t,z,\Pi)\ (I\in \cI_\la)\ \}$ 
of $\cM^{(n)}_\la(z,\Pi_\la)$. 
From Proposition \ref{triangular}, we obtain 
\begin{prop} 
 $\ds{{\rm dim}_\C \cM^{+(n)}_\la(z^{(n)},\Pi_\la)=\frac{n!}{\la_1!\cdots \la_N!}}$. 
\end{prop}



\vspace{2mm}
\noindent
{\it Remark.}\ 
For $\la=(\la_1,\cdots,\la_N)\in \N^N$,  let $x={}^t(x^{(1)}_1,\cdots,x^{(1)}_{\la^{(1)}},\cdots,x^{(N-1)}_1,\cdots,x^{(N-1)}_{\la^{(N-1)}})\in \C^M$. 
From Proposition \ref{quasiperiod} one can deduce a symmetric integral $M\times M$ matrix $N$ and a vector $\xi\in (\C/r\Z)^M$, 
which yield  
 the following quadratic form $N(x)={}^tx N x$ and the linear form $\xi(x)={}^tx \xi$. 
\be
&&N(x)=\sum_{l=1}^{N-2}\sum_{a=1}^{\la^{(l)}}\sum_{b=1}^{\la^{(l+1)}}(x^{(l)}_a-x^{(l+1)}_b)^2+n\sum_{a=1}^{\la^(N-1)}(x^{(N-1)}_a)^2,\\
&&\xi(x)=-\sum_{l=1}^{N-1}\sum_{a=1}^{\la^{(l)}}x^{(l)}_a((P+h)_{l,l+1}+\la_{l+1})-\sum_{a=1}^{\la^{(N-1)}}\sum_{k=1}^nx^{(N-1)}_au_k.
\en
Then by Appel-Humbert theorem\cite{LB}, a pair $(N,\xi)$  characterizes a line bundle $\cL(N,\xi)\ :(\C^M\times \C)/\Lambda^M\to \C^M$, where 
$\Lambda=r\Z+r\Z \tau$, with action 
\be
\omega\cdot(x,\eta)=(x+\omega,e_{\omega}(x)\eta),\qquad \omega\in \Lambda^M,\ x\in \C^M,\ \eta\in \C,
\en
and cocycle 
\be
e_{nr+mr\tau}(x)=(-1)^{{}^tnNn}(-e^{i\pi \tau})^{{}^tmNm}e^{\frac{2\pi i}{r}{}^tm(Nx+\xi)},\qquad n,m\in \Z^M.
\en 
Moreover $\Theta^+_\la(z,\Pi_\la):={\rm Span}_\C\{\ W_I(t,z,\Pi)\ (I\in \cI_\la)\ \}$ is a space of sections of $\cL(N,\xi)$. 
See for example, \cite{FRV}.

\subsection{Shuffle algebra structure}\lb{shufflealgstr}
Consider a graded $\C$-vector space
\be
&&\cM(z,\Pi)=
\bigoplus_{n\in\N}\bigoplus_{\la\in \N^N\atop |\la|=n}\cM^{(n)}_\la(z^{(n)},\Pi_\la)
\en
with $\cM^{(0)}_{(0,\cdots,0)}(z^{(0)},\Pi)=\C 1$.  

\begin{dfn}\lb{defstarprod}
For $F(t;z^{(m)}, \Pi_I)\in \cM^{(m)}_\la(z^{(m)},\Pi_\la)$,  $G(t';{z'}^{(n)},\Pi'_{I'})\in {\cM^{(n)}}_{\la'}({z'}^{(n)},\Pi'_{\la'})$, 
we define the bilinear product $\star$ on $\cM(z,\Pi)$ by 
\bea
&&(F\star G)(t^{(1)}_1,\cdots,  t^{(1)}_{\la^{(1)}+{\la'}^{(1)}}, \cdots, t^{(N-1)}_1,\cdots,  t^{(N-1)}_{\la^{(N-1)}+{\la'}^{(N-1)}};  z_1,\cdots,z_{m+n},
{\Pi}_{I+I'}
)\nn\\
&&:=\frac{1}{\prod_{l=1}^{N-1}\la^{(l)}!\la^{'(l)}!}{\rm Sym}^{(1)}\cdots {\rm Sym}^{(N-1)} \left[
F(t,z,\Pi_I q^{-2\sum_{j=1}^n<\bep_{\mu_j'},h>})\ G(t',z',\Pi'_{I'})\ {\Xi}(t,t',z,z')
\right], \nn\\
&&\lb{shufflprod}
\ena
where $I'=I'_{\mu_1'\cdots\mu'_n}$ and 
\be
&&{\Xi}(t,t',z,z')=\prod_{l=1}^{N-1}
\prod_{a=1}^{\la^{(l)}}\left(\prod_{b=1}^{\la^{'(l+1)}}\frac{[{v'_b}^{(l+1)}-v^{(l)}_a]}{[{v'_b}^{(l+1)}-v^{(l)}_a+1]}
\prod_{c=1}^{\la^{'(l)}}\frac{[{v'_c}^{(l)}-v^{(l)}_a+1]}{[{v'_c}^{(l)}-v^{(l)}_a]}\right). 
\en
In the LHS of \eqref{shufflprod},  we set $t^{(l)}_{\la^{(1)}+a}:={t'}^{(l)}_a\ (a=1,\cdots, {\la'}^{(l)})$, $z_{m+k}:=z'_{k}\ (k=1,\cdots,n)$ and 
 $\Pi_{I+I'}=\{\Pi_{\mu_k,j}\ (k=1,\cdots,m+n, j=\mu_k+1,\cdots,N) \}$, where 
 ${\Pi}_{\mu_{m+k},j}:=\Pi'_{\mu'_k,j}$ $(k=1,\cdots,n, j=\mu'_k+1,\cdots,N )$. 
\end{dfn}
This endows $\cM(z,\Pi)$ with a structure of an associative unital algebra with the unit $1$. 
In \cite{FRV}, a $\slt$ version of the $\star$-product is given. 

Let us consider the subspace of $\cM(z,\Pi)$. 
\be
&&\cM^+(z,\Pi)=\bigoplus_{n\in \N}\bigoplus_{\la\in \N^N\atop |\la|=n } \cM^{+(n)}_\la(z^{(n)},\Pi_\la). 
\en 
From \eqref{entireWmod}-\eqref{entireUmod} it turns out that  all the elements in ${\cM}^+(z,\Pi)$ satisfy the following pole and wheel conditions. 
For $F(t;z,\Pi)\in \cM^{+(n)}_\la(z^{(n)},\Pi_\la)$, 
\begin{itemize}
\item[1)] there exists an entire function $f(t;z,\Pi)\in \Theta^+_\la(z^{(n)},\Pi_\la)$ such that 
\be
&& F(t;z,\Pi)=\frac{f(t;z,\Pi)}{H_\la(t,z)}
\en
where $H_\la(t,z)$ is given in \eqref{Hfuncmod}. 
\item[2)] $f(t;z,\Pi)=0$ once $t^{(l)}_a/t^{(l+\vep)}_c=q^{2\vep}$ and $t^{(l+\vep)}_c/t^{(l)}_b=1$ for some 
$l,\vep, a,b,c$, where $\vep\in\{\pm 1\}$, $l=1,\cdots,N$, $a,b=1,\cdots,\la^{(l)}$, $c=1,\cdots,\la^{(l+\vep)}$ and $t^{(N)}_a=z_a$.

\end{itemize}

\begin{prop}\lb{shuffleMp}
The subspace $\cM^+(z,\Pi)\subset \cM(z,\Pi)$ is $\star$-closed.  
\end{prop}
A proof is given in Appendix \ref{shufflealg}.

\vspace{2mm}
\noindent
{\it Remark.}\ The subalgebra $(\cM^+(z,\Pi),\star)$ is an $A_{N-1}$-type elliptic and dynamical analogue of the shuffle algebra 
studied in \cite{Negut, FT}, where the $A^{(1)}_{N-1}$-type was discussed. $(\cM^+(z,\Pi),\star)$ is also similar to the 
elliptic algebra  discussed in \cite{FJMOP}. We will discuss the implications of the algebra $(\cM^+(z,\Pi),\star)$ in a separate paper.

\section{Elliptic $q$-KZ Equation }
In this section we consider traces of the $n$-point operators $\phi_{\mu_1\cdots\mu_n}(z_1,\cdots,z_n)$ in Sec.\ref{EWF} and show that they 
satisfy the face type elliptic $q$-KZ equation derived in \cite{Felder2,FJMMN}. Evaluating the traces explicitly we obtain  formal elliptic hypergeometric integral solutions 
to the equation. 

\subsection{Trace of the vertex  operators}
Let us consider the following trace  of the $n$-point operator.  
\bea
F^{a}(z_1,\cdots,z_n;\Pi)&=&\tr_{\F_{a,\nu}}(q^{-\kappa \hd}\ \Phi(z_1)\cdots \Phi(z_n))\nn\\
&=&\sum_{\mu_1,\cdots,\mu_n}v_{\mu_n}\tot \cdots \tot v_{\mu_1} F^{a}_{\mu_1\cdots\mu_n}(z_1,\cdots,z_n;\Pi)
\\
F^{a}_{\mu_1\cdots\mu_n}(z_1,\cdots,z_n;\Pi)&=&\tr_{\F_{a,\nu}}(q^{-\kappa \hd}\ \Phi_{\mu_1}(z_1)\cdots \Phi_{\mu_n}(z_n))
\ena
where $ \F_{a,\nu}= \F_{a,\nu}(\xi,\eta)\ (a=0,1,\cdots,N-1)$ are given in \eqref{repsp}.  To make the trace well-defined, the total weight of the operators 
in side of the trace must be zero
\be 
&&\sum_{j=1}^n\bep_{\mu_j}=0. 
\en
Since $\sum_{j=1}^N\bep_j=0$, this condition is equivalent to $\la_1=\la_2=\cdots=\la_N$
 for $I=I_{\mu_1\cdots\mu_n}\in \cI_\la$.

\begin{lem}
{
\be
&(1)&F^a_{\mu_1\mu_2\cdots\mu_n}(z_1,z_2,\cdots,q^{\kappa}z_n;\Pi)=F^{a'}_{\mu_n\mu_1\cdots\mu_{n-1}}(z_n, z_1,\cdots,z_{n-1};\Pi q^{-2<\bep_{\mu_n},h>}),\\
&(2)&F^a_{\mu_1\cdots \mu_i\mu_{i+1}\cdots\mu_n}(\cdots,z_i,z_{i+1},\cdots;\Pi)\\
&&\qquad=\sum_{\mu'_i\mu'_{i+1}} R\left(\frac{z_{i+1}}{z_i},\Pi q^{2\sum_{j=1}^{i-1}<\bep_{\mu_j},h>}\right)_{\mu_{i+1}\mu_i}^{\mu'_{i+1}\mu'_i}F^a_{\mu_1\cdots\mu'_{i+1}\mu'_{i}\cdots\mu_n}(\cdots,z_{i+1},z_{i},\cdots;\Pi),
\en
where $a'$ denotes the cyclic permutation of $a$ by $(0\ N-1\ N-2 \cdots 2\ 1)\in \gS_N$.
}
\end{lem}
\noindent
{\it Proof. }\  (1) follows from the  cyclic property of trace, $\Pi \Phi_{\mu_j}(z)=\Phi_{\mu_j}(z)\Pi q^{-2<\bep_{\mu},h>}$ , \ 
 the zero weight condition $\sum_{j=1}^n \bep_{\mu_j}=0$,   $\Phi_{\mu}(q^{\kappa }z)= q^{-\kappa \hd}\Phi_{\mu}(z)q^{\kappa \hd}$ and  
$\Phi(z) : \F_a(\xi,\eta) \to  \hV_z \tot \F_{a'}(\xi,\eta)$.\\[2mm]
 (2) follows from  the exchange relation of the vertex operators \eqref{typeIcr}. 
 \qed

\vspace{3mm}
By using the properties (1) and (2), we obtain the following statement. 
\begin{thm}
 $F^a_{\mu_1\cdots\mu_n}(z_1,\cdots,z_n;\Pi)$ satisfies the face type elliptic $q$-KZ equation
 \be
 &&\hspace{-1cm}F^a(z_1,\cdots,q^{\kappa}z_i,\cdots,z_n;\Pi)\\
&&\hspace{-0.5cm}=R^{(i+1i)}\left(\frac{q^{-\kappa}z_{i+1}}{z_i},\Pi q^{2\sum_{k=1}^{i-1}h^{{(k)}}}\right)\cdots R^{(ni)}\left(\frac{q^{-\kappa}z_{n}}{z_i},\Pi q^{2\sum_{k=1\atop \not=i}^{n-1}h^{{(k)}}}\right)\\
&&\times\Gamma_i\ R^{(1i)}\left(\frac{z_{1}}{z_i},\Pi \right)\cdots R^{(i-1i)}\left(\frac{z_{i-1}}{z_i},\Pi q^{2\sum_{k=1}^{i-2}{h^{(k)}}}\right) F^{a'}(z_1,\cdots,z_i,\cdots,z_n;\Pi). 
\en
Here $\Gamma_i$ denotes the shift operator 
\be
&&\Gamma_i f(\ \cdot\ ;\Pi)=f(\ \cdot\ ;\Pi q^{-2<\bep_\mu,h>})
\en
if $q^{2h^{(i)}} f(\ \cdot\ ;\Pi)=q^{2<\bep_\mu,h>} f(\ \cdot\ ;\Pi)$ for $f(\ \cdot\ ;\Pi)=\sum v_{\mu_n}\tot \cdots\tot v_{\mu_1}\tot f_{\mu_1\cdots\mu_n}(\ \cdot\ ;\Pi)$.    
\end{thm}

\subsection{Evaluation of the trace}
We now apply the formula in Theorem \ref{prodPhi} to the $n$-point operator $\Phi_{\mu_1}(z_1)\cdots \Phi_{\mu_n}(z_n)$, and 
 evaluate the trace $\Phi(t,z):=\tr_{\F_{a,\nu}}\left(q^{-\kappa d}\ \widetilde{\Phi}(t,z)\right)$ explicitly.
We  obtain the following result. 
\begin{thm} For $|p|<|z_1|, \cdots, |z_n|<1 $, we obtain
\bea
&&\hspace{-2cm}
F^a_{\mu_1\cdots\mu_n}(z_1,\cdots,z_n;\Pi)
=\oint_{\T^M} \underline{dt}\  \Phi(t,z)\ \omega_{\mu_1\cdots\mu_n}(t,z,\Pi),\lb{tracenpt}\\
\Phi(t,z)
&=& C_n 
(-)^{n h_{\bep_N}}\prod_{l=1}^{N-1}q^{(N-l)\la^{(l)}\frac{(P+h)_{\al_l}-1}{r}}\nn\\
&& \times\prod_{k=1}^n\left(z_k^{-\la^{(N-1)}-h_{\bep_N}+\frac{1}{r}( (P+h)_{\bep_N}+  \la^{(N-1)})}
\prod_{1\leq k\not= l\leq n}\frac{\Gamma(q^2z_k/z_l;p,q^{\kappa},q^{2N})}{\Gamma(q^{2N}z_k/z_l;p,q^{\kappa},q^{2N})}\right)
\nn\\
&&\times \prod_{l=1}^{N-1}\left[
\prod_{a=1}^{\la^{(l)}}\left( (t^{(l)}_a)^{\la_l-h_{\al_l}+\frac{1}{r}( (P+h)_{\al_l} -\la^{(l)})}
\prod_{b=1}^{\la^{(l+1)}}
\frac{\Gamma(qt^{(l)}_a/t^{(l+1)}_b;p,q^{\kappa})}{\Gamma(p^*qt^{(l)}_a/t^{(l+1)}_b;p,q^{\kappa})}\right)\right.
\nn\\
&&\left.\qquad\qquad\qquad\qquad \times 
\prod_{1\leq a<b\leq \la^{(l)}}\frac{\Gamma(p^*t^{(l)}_a/t^{(l)}_b,p^*t^{(l)}_b/t^{(l)}_a;p,q^{\kappa})}{\Gamma(t^{(l)}_a/t^{(l)}_b,t^{(l)}_b/t^{(l)}_a;p,q^{\kappa})}\right],
\ena
where
\be
C_n&=&(-)^{n\la^{_(N-1)}}(-q^{\frac{r-1}{r}})^{\frac{1}{2}\sum_{l=1}^{N-1}\la^{(l)}(\la^{(l)}-1)}\\
&&\times \left(\frac{(p;p)_\infty}{(q^2;p)_\infty}\frac{\Gamma(p;p,q^\kappa)}{\Gamma(q^2;p,q^\kappa)}\right)^{\sum_{l=1}^{N-1}\la^{(l)}}
\left(\frac{(q^{2N};p,q^{2N})_\infty}{(q^2;p,q^{2N})_\infty}\frac{\Gamma(q^2;p,q^\kappa,q^{2N})}{\Gamma(q^{2N};p,q^\kappa,q^{2N})}\right)^n.
\en
\end{thm}

\noindent
{\it Proof.}\ 
In $\tPhi(t,z)$, taking  normal ordering further among $\Phi_N$'s and $F_{N-1}$'s as well as among $F_{l+1}$'s and $F_l$'s $(l=1,\cdots,N-1)$, we obtain
\be
\Phi(t,z)&=&
:\Phi_N(z_1)
\cdots 
\Phi_N(z_n)
F_{N-1}(t^{(N-1)}_{1})\cdots F_{N-1}(t^{(N-1)}_{\la^{(N-1)}})\cdots  F_1(t_1^{(1)})\cdots F_{1}(t^{(1)}_{\la^{(1)}}):
\nn\\
&&\times \prod_{1\leq l<m\leq n}<\Phi_N(z_l)\Phi_N(z_m)>^{Sym}
\prod_{l=1}^{N-1}\prod_{1\leq a<b\leq \la^{(l)}}<F_l(t^{(l)}_a)F_l(t^{(l)}_b)>^{Sym}\\
&&{\times\prod_{l=1}^n\prod_{a=1}^{\la^{(N-1)}}<\Phi_N(z_l)F_{N-1}(t^{(N-1)}_a)>\prod_{l=1}^{N-2}\prod_{a=1}^{\la^{(l+1)}}\prod_{b=1}^{\la^{(l)}}
<F_{l+1}(t^{(l+1)}_a)F_l(t^{(l)}_b)>},
\en
where
\be
<\Phi_N(z)F_{N-1}(t)>&=&-z^{-(r-1)/r}\frac{(pq^{-1}t/z;p)_\infty}{(qt/z;p)_\infty},\\
<F_{l+1}(t_a)F_l(t_b)>&=&t_a^{-(r-1)/r}\frac{(pq^{-1}t_b/t_a;p)_\infty}{(qt_b/t_a;p)_\infty} .
\en
Then using the formulas in Theorem \ref{bareboson} and \ref{levelone}, the  trace of the normal ordered operator $:\Phi_N(z_1)\cdots  F_{1}(t^{(1)}_{\la^{(1)}}):$  in $\tPhi(t,z)$ 
can be evaluated, for example,  by  the coherent state method. Combining the result with all the OPE coefficients in $\tPhi(t,z)$ 
we obtain the desired result. 
\qed 

\vspace{2mm}
\noindent
{\it Remark.} The integrand of \eqref{tracenpt} is a single valued function of $t^{(1)}_1,\cdots,t^{(N-1)}_{\la^{(N-1)}}$ due to Proposition \ref{meroW}. 

\vspace{2mm}
\noindent
{\it Remark.} As discussed in \cite{TV97} for the trigonometric $\slth$ case, one can  specify the cycles of the integral \eqref{tracenpt}
by further inserting the other elliptic weight functions 
$\omega^{\kappa}_{I}(t,z,\Upsilon)$, which are the same as $\omega_I(t,z,\Pi)$ in \eqref{def:omega} except for replacing $p$ by $q^{\kappa}$,
i.e. replacing all the theta functions with elliptic nome $p$ by those with $q^\kappa$, and 
$\Pi_{j,k}$ by the other dynamical parameters $\Upsilon_{j,k}$. 
Then we have the following elliptic hypergeometric paring for $I,J\in \cI_\la$ 
\bea
&&I(\omega^\kappa_I,\omega_J)=\oint_{\T^M} \underline{dt}\  \Phi(t,z)\ \omega^\kappa_I(t,z,\Upsilon)\omega_{J}(t,z,\Pi). 
\ena
We will discuss the details of  such integral in elsewhere. 

\vspace{2mm}
\noindent
{\it Remark.}\ The elliptic algebra $U_{q,p}(\slnh)$ at the level 1 is deeply related to the deformed $W$-algebra $\cW_{q,t}(\sln)$ in \cite{FF,AKOS} 
: identifying $p$,  $p^*=p q^{-2}$ in $U_{q,p}(\slnh)$ with $q$, $t$ in $\cW_{q,t}(\sln)$, respectively,   
the elliptic currents $F_j(z), E_j(z)\ (j=1,\cdots,N-1)$ can be identified with the screening currents  $S^+_j(z), S^-_j(z)$, respectively\cite{K98,JKOS},
 the bare type I vertex 
$\Phi_N(z)$ and the type II vertex $\Psi^*_N(z)$ with certain deformed primary fields in $\cW_{q,t}(\sln)$\cite{Konno1314}
 as well as the space $\F_{a,\nu}$ with the Verma module of $\cW_{q,t}(\sln)$\cite{FKO}. 
 Hence $\Phi(t,z)$ can be regarded as an (elliptic) correlation function of  $\cW_{q,t}(\sln)$. 
 Recently, Aganagic, Frenkel and Okounkov\cite{AFO} proposed 
  the notion of quantum $q$-Langlands duality, which states the correspondence between the solutions to the trigonometric $q$-KZ equation 
associated with $U_q(\gh)$ and the correlation functions of $\cW_{q,t}({}^L\g)$ (${}^L\g$ is the Langlands dual Lie algebra of $\g$). 
If one identifies  the elliptic weight functions $\omega_{\mu_1\cdots\mu_n}(t,z,\Pi)$ with the stable envelopes, 
the formula \eqref{tracenpt} seems to give an elliptic analogue of the formulas, for example, (5.11) in \cite{AFO}. In this sense 
the quantum $q$-Langlands duality seems to become more transparent  in the elliptic algebra level. We will discuss this point in a separate paper.

\section*{Acknowledgements}
The author would like to thank Vassily Gorbounov, Masahiko Ito, Michio Jimbo, Masatoshi Noumi, Andrei Okounkov and Vitaly Tarasov 
for discussions. He is also grateful to 
Igor Frenkel, Giovanni Felder, Omar Foda and Alexander Varchenko for  useful conversations and their interests. 
His research is supported in part by the Grant-in -Aid for Scientific Research (C) 17K05195, JSPS.

\begin{appendix}

\section{Elliptic Algebra $U_{q,p}(\slnh)$}\lb{secuqpBn}
In this appendix,  we summarize some basic facts on the elliptic algebra $U_{q,p}(\slnh)$. 

\subsection{Definition}
In this subsection, let $q=e^{\hbar}\in \C[[\hbar]]$ and $p$ be an indeterminate. 

\begin{dfn}\cite{FKO}\lb{defUqp}
The elliptic algebra $U_{q,p}(\slnh)$ is a topological algebra over $\FF[[p]]$ generated by 
$e_{j,m}, f_{j,m}, \al_{j,n}, K^\pm_j$,  
$(1\leq j\leq N-1, m\in \Z, n\in \Z_{\not=0})$, ${d}$ and the central element $c$. 
We assume $K^\pm_{j}$ are invertible and set 
\be
&&e_j(z)=\sum_{m\in \Z} e_{j,m}z^{-m}
,\quad f_j(z)=\sum_{m\in \Z} f_{j,m}z^{-m},\lb{deffn}\\
&&{\psi}_j^+(
q^{-\frac{c}{2}}
z)=K^+_{j}\exp\left(-(q-q^{-1})\sum_{n>0}\frac{\al_{j,-n}}{1-p^n}z^n\right)
\exp\left((q-q^{-1})\sum_{n>0}\frac{p^n\al_{j,n}}{1-p^n}z^{-n}\right),\\
&&{\psi}_j^-(q^{\frac{c}{2}}z)=K^-_{j} \exp\left(-(q-q^{-1})\sum_{n>0}\frac{p^n\al_{j,-n}}{1-p^n}z^n\right)
\exp\left((q-q^{-1})\sum_{n>0}\frac{\al_{j,n}}{1-p^n}z^{-n}\right). 
\en
We call $e_j(z), f_j(z), \psi^\pm_j(z)$ the elliptic currents. 
The  defining relations are as follows. For $g(P), g(P+h)\in \cM_{H^*}$, 
\bea
&&g({P+h})e_j(z)=e_j(z)g({P+h}),\quad g({P})e_j(z)=e_j(z)g(P-<Q_{\al_j},P>),\lb{ge}\\
&&g({P+h})f_j(z)=f_j(z)g(P+h-<{\al_j},P+h>),\quad g({P})f_j(z)=f_j(z)g(P),\lb{gf}
\\
&&[g(P), \al_{i,m}]=[g(P+h),\al_{i,n}]=0,\lb{gboson}\\ 
&&g({P})K^\pm_j=K^\pm_jg(P-<Q_{\al_j},P>),\\
&&\ g({P+h})K^\pm_j=K^\pm_jg(P+h-<Q_{\al_j},P>),\lb{gKpm}
\\
&&[d, g(P+h)]=[d, g(P)]=0,
\quad\lb{dg}\\
&& [d, \al_{j,n}]=n\al_{j,n},\quad [d, e_j(z)]=-z\frac{\partial}{\partial z}e_j(z), \quad  [d, f_j(z)]=-z\frac{\partial}{\partial z}f_j(z),\quad \lb{dedf}\\
&&K_i^{\pm}e_j(z)=q^{\mp a_{ij}}e_j(z)K_i^{\pm},\quad 
K_i^{\pm}f_j(z)=q^{\pm a_{ij}}f_j(z)K_i^{\pm},
\\
&&[\al_{i,m},\al_{j,n}]=\delta_{m+n,0}\frac{[a_{ij}m]_q
[cm]_q
}{m}
\frac{1-p^m}{1-p^{*m}}
q^{-cm}
,\lb{ellboson}\\
&&
[\al_{i,m},e_j(z)]=\frac{[a_{ij}m]_q}{m}\frac{1-p^m}{1-p^{*m}}
q^{-cm}z^m e_j(z),
\lb{bosonve}\\
&&
[\al_{i,m},f_j(z)]=-\frac{[a_{ij}m]_q}{m}z^m f_j(z)
,\lb{bosonvf}
\ena
\bea
&&
z_1 \frac{(q^{a_{ij}}z_2/z_1;p^*)_\infty}{(p^*q^{-a_{ij}}z_2/z_1;p^*)_\infty}e_i(z_1)e_j(z_2)=
-z_2 \frac{(q^{a_{ij}}z_1/z_2;p^*)_\infty}{(p^*q^{-a_{ij}}z_1/z_2;p^*)_\infty}e_j(z_2)e_i(z_1),\lb{ee}\\
&&
z_1 \frac{(q^{-a_{ij}}z_2/z_1;p)_\infty}{(pq^{a_{ij}}z_2/z_1;p)_\infty}f_i(z_1)f_j(z_2)=
-z_2 \frac{(q^{-a_{ij}}z_1/z_2;p)_\infty}{(pq^{a_{ij}}z_1/z_2;p)_\infty}f_j(z_2)f_i(z_1),\lb{ff}\\
&&[e_i(z_1),f_j(z_2)]=\frac{\delta_{i,j}}{q-q^{-1}}
\left(\delta(
q^{-c}
z_1/z_2)
\psi_j^-(
q^{\frac{c}{2}}
z_2)-
\delta(
q^c
z_1/z_2)
\psi_j^+(
q^{-\frac{c}{2}}
z_2)
\right),\lb{eifj}
\ena
\bea
&&\frac{(p^*q^2z_2/z_1;p^*)_\infty}
{(p^*q^{-2}z_2/z_1;p^*)_\infty
} \left\{
\frac{(p^*q^{-1}z/z_1;p^*)_\infty 
(p^*q^{-1}z/z_2;p^*)_\infty}
{(p^*qz/z_1;p^*)_\infty 
(p^*qz/z_2;p^*)_\infty}e_i(z_1)
e_i(z_2)
e_j(z)\right.\nonumber\\
&&\qquad\qquad -\left.[2]_q\frac{(p^*q^{-1}z/z_1;p^*)_\infty 
(p^*q^{-1}z_2/z;p^*)_\infty}
{(p^*qz/z_1;p^*)_\infty 
(p^*qz_2/z;p^*)_\infty}
e_i(z_{1})e_j(z)e_i(z_{2})
\right.\nonumber\\
&&\qquad\qquad +\left.
\frac{(p^*q^{-1}z_1/z;p^*)_\infty 
(p^*q^{-1}z_2/z;p^*)_\infty}
{(p^*qz_1/z;p^*)_\infty 
(p^*qz_2/z;p^*)_\infty}
e_j(z)e_i(z_{1})
e_i(z_{2})
\right\}+(z_1 \leftrightarrow z_2)=0,\nn\\
&&\lb{serre}\\
&&\frac{(pq^{-2}z_2/z_1;p)_\infty}
{(pq^{2}z_2/z_1;p)_\infty
} \left\{
\frac{(pq z/z_1;p)_\infty 
(pq z/z_2;p)_\infty}
{(pq^{-1} z/z_1;p)_\infty 
(pq^{-1} z/z_2;p)_\infty}
f_i(z_1)f_i(z_2)f_j(z)\right.\nonumber\\
&&\qquad\qquad
-\left.[2]_q\frac{(pq z/z_1;p)_\infty 
(pq z_2/z;p)_\infty}
{(pq^{-1} z/z_1;p)_\infty 
(pq^{-1} z_2/z;p)_\infty}
f_i(z_{1})f_j(z)f_i(z_{2})
\right.\nonumber
\\
&&\qquad\qquad+\left.
\frac{(pq z_1/z;p)_\infty 
(pq z_2/z;p)_\infty}
{(pq^{-1} z_1/z;p)_\infty 
(pq^{-1} z_2/z;p)_\infty}
f_j(z)f_i(z_{1})
f_i(z_{2})
\right\} +(z_1 \leftrightarrow z_2)=0,\nn\\
&&~~~~~~~~~~~\qquad\qquad{\rm for}~~~|i-j|=1,\lb{serrf}
\ena
where $p^*=p
q^{-2c}
$ and $\delta(z)=\sum_{n\in \Z}z^n$.  
\end{dfn}
We treat these relations 
as formal Laurent series in $z, w$ and $z_j$'s. 
All the coefficients are well-defined in the $p$-adic topology.

\subsection{The orthonormal basis type elliptic bosons}

Let us define the orthonormal basis type elliptic bosons $\cE^{ j}_m \ (1\leq j\leq N, m\in \Z_{\not=0})$ \cite{FKO} by
\bea
&&\cE^{ j}_m= q^{ jm}\frac{C_m}{q-q^{-1}}\left(-q^{- Nm}\sum_{k=1}^{j-1}[km]_q\al_{k,m}
+\sum_{k=j}^{N-1} [ (N-k ) m]_q\al_{k,m} \right) \lb{orthoboson}
\ena
with
\be
&&C_m=\frac{1}{[m]_q^2[N m]_q}.
\en
One can show the following relations. 
\begin{prop}\lb{cEcE}
\bea
&&\sum_{j=1}^Nq^{(j-1)m}\cE^{ j}_m=0,\\
&&[\cE^{ j}_m,\cE^{ j}_n]
=\delta_{m+n,0}\frac{[cm]_q[\eta m]_q [(N-1)m]_q}{m(q-q^{-1})^2[m]_q^3[Nm]_q }\frac{1-p^m}{1-p^{*m}}q^{-cm}, 
\\
&&[\cE^{ j}_m,\cE^{ k}_n]
=-\delta_{m+n,0}q^{({\rm sgn}(k-j) N -k+j)m}
\frac{[cm]_q}{m(q^m-q^{-m})^2[N m]_q}\frac{1-p^m}{1-p^{*m}}q^{-cm},
\ena
where
\be
{\rm sgn}(l-j)=\left\{\mmatrix{+&(l>j)\cr
0&(l=j)\cr 
                              -&(l<j). \cr }\right.
\en
\end{prop}

Note also
\begin{prop}\lb{alphaAA}
\noindent 
\bea
&&\al_{j,m}=[m]_q^2(q-q^{-1})(\cE^{ j}_m-q^{- m}\cE^{j+1}_m).
\ena
\end{prop}
Furthermore if we set
\be
&&A^j_m=(q^m-q^{-m})\sum_{k=1}^jq^{(k-j-1)m}\cE^{ k}_m,
\en
then we have
\be
&&[\al_{i,m},A^j_m]=-\delta_{i,j}\delta_{m+n,0}\frac{[cm]_q}{m}\frac{1-p^m}{1-p^{*m}}q^{-cm} \quad (1\leq i,j\leq N-1)
\en
Hence we call $A^j_m$  the fundamental weight type bosons.

\subsection{The elliptic currents $k_{ j}(z)$}
Let us set 
\bea
&&\psi_j(z)=:\exp\left\{(q-q^{-1})\sum_{m\not=0}\frac{\al_{j,m}}{1-p^m}p^mz^{-m}\right\}:.
\ena
Here $:\quad :$ denotes the normal ordering  defined by
\be
&&:\al_{j,m}\al_{k,n}:=\left\{\mmatrix{\al_{j,m}\al_{k,n}&\mbox{if}\ m\leq n\cr
                                                  \al_{k,n}\al_{j,m}&\mbox{if}\ m> n\cr
}\right.
\en
for $1\leq j,k\leq N-1$. 
Then the elliptic currents $\psi_j^\pm(z)$ in Definition \ref{defUqp} can be written as 
\bea
&&\psi^+_j(q^{-\frac{c}{2}}z)=K^+_j\psi_j(z),\qquad \psi^-_j(q^{-\frac{c}{2}}z)=K^-_j\psi_j(pq^{-c}z).
\ena
Let us introduce the new elliptic currents $k_{  j}(z) \,  (j\in I \cup\{N\})$ associated with $\cE^{ j}_m$ by
\bea
k_{  j}(z) &=& :\exp\left\{\sum_{m\not=0} \frac{[m]_q^2(q-q^{-1})^2}{1-p^m}p^m\cE^{ j}_mz^{-m} \right\}:
\ena
Then from Proposition \ref{alphaAA} we obtain the following relations.  
\begin{prop}\lb{psikk}
\bea
\psi_j(z)&=&\varrho
k_{j}(z)k_{j+1}(qz)^{-1}, 
\ena
where
\bea
&&
\varrho=\frac{(p;p)_\infty(p^*q^2;p^*)_\infty}{(p^*;p^*)_\infty(pq^2;p)_\infty}.\lb{defkappa}
\ena
\end{prop}

In addition, from Proposition \ref{cEcE} we obtain the following commutation relations. 
\begin{thm}\lb{kk}
\be
&&k_{  j}(z_1)k_{ j}(z_2)=\frac{\tilde{\rho}^{+*}(z)}{\tilde{\rho}^+(z)}k_{ j}(z_2)k_{ j}(z_1),
\qquad  (1\leq j\leq N),\\
&&k_{ j}(q^{ j}z_1)k_{ k}(q^{ k}z_2)=\frac{\tilde{\rho}^{+*}(z)}{\tilde{\rho}^{+}(z)}\frac{\Theta_{p^*}(q^{-2}z)\Theta_{p}(z)}{\Theta_{p^*}(z)\Theta_{p}(q^{-2}z)} k_{ k}(q^{ k}z_2)k_{j}(q^{ j}z_1) \qquad (1\leq j<k\leq N),
\en
where $z=z_1/z_2$,  $\tilde\rho^+(z)$ is given in \eqref{defrhotilde} and $\rho^{+*}(z)=\rho^+(z)|_{p\mapsto p^*}$.

\end{thm}

\begin{prop}\lb{kef}
\be
&&k_{j}(z_1)e_j(z_2)=\frac{\Theta_{p^*}(q^{-c}z)}{\Theta_{p^*}(q^{-c- 2}z)}
e_j(z_2)k_{  j}(z_1)\qquad (1\leq j\leq N),\\
&&k_{ j}(z_1)e_{j-1}(z_2)=\frac{\Theta_{p^*}(q^{-c-1}z)}{\Theta_{p^*}(q^{-c+ 1}z)}
e_{j-1}(z_2)k_{ j}(z_1)\qquad (2\leq j\leq N),\\
&&k_{ j}(z_1)e_k(z_2)=e_k(z_2)k_{ j}(z_1)\qquad(k\not=j,j-1), \\
&&k_{ j}(z_1)f_j(z_2)=\frac{\Theta_{p}(q^{-2}z)}{\Theta_{p}(z)}
f_j(z_2)k_{ j}(z_1)\qquad (1\leq j\leq N),\\
&&k_{ j}(z_1)f_{j-1}(z_2)=\frac{\Theta_{p}(qz)}{\Theta_{p}(q^{-1}z)}
f_{j-1}(z_2)k_{ j}(z_1)\qquad (2\leq j\leq N),\\
&&k_{ j}(z_1)f_k(z_2)=f_k(z_2)k_{ j}(z_1)\qquad(k\not=j,j-1).
\en
\end{prop}


\subsection{Modified elliptic currents}\lb{MEC}

The $R$ matrix \eqref{ellR} is gauge equivalent to Jimbo-Miwa-Okado's $A_{N-1}^{(1)}$ face type Boltzmann weight\cite{JMO} 
 and  conveniently  
expressed by using Jacobi's theta function \eqref{thetafull}.
However one drawback is that  Jacobi's theta is  accompanied by 
the fractional power of $z$. 
In order to deal with this one needs to introduce the following modifications of the elliptic currents\cite{KK03}.

\begin{dfn}\lb{extragen}
We 
introduce the new generators $e^{\pm \zeta_{\bep_j}}\ (1\leq j\leq N)$ satisfying 
\bea
&&
e^{Q_{\bar{\epsilon}_{j}}}e^{ Q_{\bar{\epsilon}_{k}}}=
q^{\left(\frac{1}{r}-\frac{1}{r^*}\right){\rm sgn}(j-k)} e^{Q_{\bar{\epsilon}_{k}}}e^{ Q_{\bar{\epsilon}_{j}}},\lb{central2}\\
&&e^{Q_{\bar{\epsilon}_{j}}}e^{\zeta_{\bep_k}}=
q^{\frac{1}{r}{\rm sgn}(j-k)}e^{\zeta_{\bep_k}}e^{Q_{\bar{\epsilon}_{j}}},\lb{central3}\\
&&e^{\zeta_{\bep_j}}e^{\zeta_{\bep_k}}=q^{\frac{1}{r}{\rm sgn}(j-k)}e^{\zeta_{\bep_k}}e^{\zeta_{\bep_j}},\lb{central4}\\
&&[P_{\bar{\epsilon}_j}, e^{\pm \zeta_{\bep_k}}]=0,\qquad e^{\pm\sum_{j=1}^N \zeta_{\bep_j}}=1,\lb{central5}\\
&&
[e^{\pm \zeta_{\bep_j}}, \C[\cQ]]=
[e^{\pm \zeta_{\bep_j}}, U_{q,p}(\slnh)]=0.
\lb{central6}
\ena
\end{dfn}
Let us set $\zeta_j=\zeta_{\bep_j}-\zeta_{\bep_{j+1}}$, we have 
\begin{prop}
\bea
&&e^{Q_{\al_{j}}}e^{ Q_{\al_{k}}}=
q^{\left(\frac{1}{r}-\frac{1}{r^*}\right)\left(\delta_{j,k+1}-
\delta_{j,k-1}\right)}e^{Q_{\al_{k}}}e^{ Q_{\al_{j}}},\lb{qajqak}\\
&&e^{Q_{\bep_{j}}}e^{ Q_{\al_{k}}}=q^{
\left(\frac{1}{r}-\frac{1}{r^*}\right)\left(\delta_{j,k}+
\delta_{j,k+1}\right)}e^{ Q_{\al_{k}}} e^{Q_{\bep_{j}}},\lb{qbejqak}\\
&&e^{Q_{\bar{\epsilon}_{j}}}e^{\zeta_k}=
q^{\frac{1}{r}\left(\delta_{j,k}+\delta_{j,k+1}\right)}e^{\zeta_k}e^{Q_{\bar{\epsilon}_{j}}},\lb{qbejbak}\\
&&e^{Q_{\al_{j}}}e^{\zeta_k}=
q^{\frac{1}{r}\left(\delta_{j,k+1}-\delta_{j,k-1}\right)}e^{\zeta_k}e^{Q_{\al_{j}}},\lb{qajbak}\\
&&e^{\zeta_j}e^{ \zeta_k}=q^{\frac{1}{r}\left(\delta_{j,k+1}-
\delta_{j,k-1}\right)}e^{\zeta_k}e^{ \zeta_j},
\lb{bajbak}\\
&&[P_{\bar{\epsilon}_k},e^{\pm\zeta_j} ]=[e^{\pm \zeta_j}, U_{q,p}(\slnh)]=0.\lb{pbejbak}
\ena
\end{prop}
\begin{dfn}\lb{modellcurrents}
We  define the modified elliptic currents as follows.
\be
\!\!\!&&E_{ j}(z)=e_j(z)e^{\zeta_j} (q^{N-j}z)^{-\frac{ P_{\al_j}-1}{r^*}}\qquad (1\leq j\leq N-1),\\
\!\!\!&&F_{ j}(z)=f_j(z)e^{-\zeta_j} (q^{N-j}z)^{\frac{ (P+ h)_{\al_j}-1}{r}}\qquad (1\leq j\leq N-1),\\
&&K^+_{j}(z)=k_{+j}(q^{j-N+1} z)e^{-Q_{\bep_j}}q^{-h_{\bep_j}}(q^{-r+1}z)^{-\frac{1}{r^*}(P_{\bep_j}-1)
-\frac{1}{r}((P+ h)_{\bep_j}-1)},\\
&&K^-_{ j}(z)=K^+_{ j}\left(pq^{-c}z\right).  
\en
for $1\leq j\leq N$. 
We also set
\bea
H^{\pm}_{j}(z)
&=&
\varrho K^\pm_{j}\left(q^{{N-j-1}}q^{\frac{c}{2}}z\right)K^\pm_{j+1}\left(q^{{N-j-1}}q^{\frac{c}{2}}z\right)^{-1}
\\
&=&\psi^{\pm}_{j}(z)(K_j^\pm)^{-1}e^{-Q_{\al_j}} q^{\mp h_j}
(q^{N-j}q^{\pm\frac{c}{2}}z)^{-\frac{1}{r^*}( P_{\al_j}-1)
+\frac{1}{r}((P+ h)_{\al_j}-1)},
\qquad (1\leq j\leq N-1).\nn
\ena
and
\bea
\hd=d+\frac{1}{2r^*}\sum_{j=1}^N(P_j+2)P^j-\frac{1}{2r}\sum_{j=1}^N((P+h)_j+2)(P+h)^j.\lb{dtilde}
\ena
\end{dfn}
Then the defining relations \eqref{dedf}--\eqref{serrf} 
of $U_{q,p}(\slnh)$ can be 
rewritten as follows in the sense of analytic continuation. 
\begin{prop}\lb{defrelfull}
\bea
&&[h_i,\al_{j,n}]=0,\quad [h_i,E_j(z)]=a_{ij} E_j(z),\quad 
[h_i,F_j(z)]=-a_{ij} F_j(z), 
\lb{u2}\\
&&[{\hd},h_i]=0,\quad [{\hd},\al_{i,n}]= n \al_{i,n},\quad\\ 
&&[{\hd},E_{i}(z)]=\left(-z\frac{\partial}{\partial z}+\frac{1}{r^*}\right)
E_i(z) , \quad
[{\hd},F_{i}(z)]=\left(-z\frac{\partial}{\partial z}+\frac{1}{r}\right)
F_i(z) ,
\lb{u3}
\\
&&[\al_{i,m},\al_{j,n}]=\delta_{m+n,0}\frac{[a_{ij}m]_q
[cm]_q
}{m}
\frac{1-p^m}{1-p^{*m}}
q^{-cm}
,\lb{ellbosonmod}\\
&&
[\al_{i,m},E_j(z)]=\frac{[a_{ij}m]_q}{m}\frac{1-p^m}{1-p^{*m}}
q^{-cm}z^m E_j(z),
\lb{bosonvE}\\
&&
[\al_{i,m},F_j(z)]=-\frac{[a_{ij}m]_q}{m}z^m F_j(z)
,\lb{bosonvF}
\ena
\bea
&& \left[u-v-\frac{a_{ij}}{2}\right]^* E_i(z)E_j(w)
=\left[u-v+\frac{a_{ij}}{2}\right]^* 
E_j(w)E_i(z),
\lb{u7}\\
&& \left[u-v+\frac{a_{ij}}{2} \right] F_i(z)F_j(v)
=\left[u-v-\frac{a_{ij}}{2}\right]
F_j(w)F_i(z),
\lb{u8}\\
&&[E_i(z),F_j(w)]=\frac{\delta_{i,j}}{q-q^{-1}}
\left(\delta\bigl(q^{-c}\frac{z}{w}\bigr)H^-_{i}(q^{c/2}w)
-\delta\bigl(q^{c}\frac{z}{w}\bigr)H^+_{i}(q^{-c/2}w)
\right),\lb{u9}
\ena
\bea
&&z_1^{-\frac{1}{r^*}}
\frac{(p^*q^2z_2/z_1;p^*)_\infty}
{(p^*q^{-2}z_2/z_1;p^*)_\infty
}\left\{
({z_2}/{z})^{\frac{1}{r^*}}
\frac{(p^*q^{-1}z/z_1;p^*)_\infty 
(p^*q^{-1}z/z_2;p^*)_\infty}
{(p^*qz/z_1;p^*)_\infty 
(p^*qz/z_2;p^*)_\infty}E_i(z_1)E_i(z_2)E_j(z)\right.\nonumber\\
&&
-\left.[2]_q\frac{(p^*q^{-1}z/z_1;p^*)_\infty 
(p^*q^{-1}z_2/z;p^*)_\infty}
{(p^*qz/z_1;p^*)_\infty 
(p^*qz_2/z;p^*)_\infty}
E_i(z_{1})E_j(z)E_i(z_{2})
\right.\nn
\\
&&+\left.
(z/z_1)^{\frac{1}{r^*}}\frac{(p^*q^{-1}z_1/z;p^*)_\infty 
(p^*q^{-1}z_2/z;p^*)_\infty}
{(p^*qz_1/z;p^*)_\infty 
(p^*qz_2/z;p^*)_\infty}
E_j(z)E_i(z_{1})E_i(z_{2})
\right\}+(z_1 \leftrightarrow z_2)=0,\nn\\
&&\label{e20}
\\
&&z_1^{\frac{1}{r}}
\frac{(pq^{-2}z_2/z_1;p)_\infty}
{(pq^{2}z_2/z_1;p)_\infty
}\left\{
(z/z_2)^{\frac{1}{r}}
\frac{(pq z/z_1;p)_\infty 
(pq z/z_2;p)_\infty}
{(pq^{-1} z/z_1;p)_\infty 
(pq^{-1} z/z_2;p)_\infty}
F_i(z_1)F_i(z_2)F_j(z)\right.\nonumber\\
&&
-\left.[2]_q\frac{(pq z/z_1;p)_\infty 
(pq z_2/z;p)_\infty}
{(pq^{-1} z/z_1;p)_\infty 
(pq^{-1} z_2/z;p)_\infty}
F_i(z_{1})F_j(z)F_i(z_{2})
\right.\nn
\\
&&+\left.
(z_1/z)^{\frac{1}{r}}
\frac{(pq z_1/z;p)_\infty 
(pq z_2/z;p)_\infty}
{(pq^{-1} z_1/z;p)_\infty 
(pq^{-1} z_2/z;p)_\infty}
F_j(z)F_i(z_{1})
F_i(z_{2})
\right\} + (z_1 \leftrightarrow z_2)=0
\quad (|i-j|=1).\nonumber\\
&&~~~~~~~~~~~~~~~~~\label{e21}
\ena

\end{prop}
In addition, one can rewrite the formulas in Theorem \ref{kk} and Proposition \ref{kef} as follows.
\begin{prop}\lb{relKK}\lb{relEK}
\be
&&K^+_{ j}(z_1)K^+_{ j}(z_2)=
\frac{{\rho}^{+*}(z_1/z_2)}{{\rho}^+(z_1/z_2)}K^+_{ j}(z_2)K^+_{ j}(z_1),\\
&&K^+_{j}(z_1)K^+_{l}(z_2)=
\frac{{\rho}^{+*}(z_1/z_2)}{{\rho}^+(z_1/z_2)}\frac{[u_1-u_2-1]^*[u_1-u_2]}{[u_1-u_2]^*[u_1-u_2-1]}K^+_{l}(z_2)K^+_{j}(z_1) \quad (1\leq j < l\leq N),\\
&&K^+_{j}(z_1)E_j(z_2)=
\frac{\left[u_1-u_2+\frac{j-N-c+1}{2}\right]^*}{\left[u_1-u_2+\frac{j-N-c+1}{2}-1\right]^*}E_j(z_2)K^+_{j}(z_1)\qquad (1\leq j\leq N),
\\
&&K^+_{j+1}(z_1)E_{j}(z_2)=
\frac{\left[u_1-u_2+\frac{j-N+1-c}{2}\right]^*}{\left[u_1-u_2+\frac{j-N+1-c}{2}+1\right]^*}E_{j}(z_2)K^+_{j+1}(z_1)\qquad (1\leq j\leq N-1),\\
&&K^+_l(z_1)E_j(z_2)=E_j(z_2)K^+_l(z_1)\qquad (l\not=j,j+1),
\\
&&K^+_{j}(z_1)F_j(z_2)=
\frac{\left[u_1-u_2+\frac{j-N+1}{2}-1\right]}{\left[u_1-u_2+\frac{j-N+1}{2}\right]}F_j(z_2)K^+_{j}(z_1)\qquad (1\leq j\leq N),\\
&&K^+_{j+1}(z_1)F_{j}(z_2)=\frac{\left[u_1-u_2+\frac{j-N+1}{2}+1\right]}{\left[u_1-u_2+\frac{j-N+1}{2}\right]}F_{j}(z_2)K^+_{j+1}(z_1)\qquad (1\leq j\leq N-1),\\
&&K^+_l(z_1)F_j(z_2)=F_j(z_2)K^+_l(z_1)\qquad (l\not=j,j+1).
\en
\end{prop}


\subsection{The half currents and the $L$-operators}
We define the half currents of $U_{q,p}(\slnh)$
as follows.
\begin{dfn}\lb{def:halfcurrents}\cite{KK03}
Let us assume $|p|<|z|<1$. We set 
\begin{eqnarray}
F_{j,l}^+(z)&=&a_{j,l}\oint_{\T^{l-j}
}
\prod_{m=j}^{l-1}\frac{dt_m}{2\pi i t_m}
F_{l-1}(t_{l-1})F_{l-2}(t_{l-2})
\cdots F_{j}(t_{j})\nonumber\\
&&\hspace{-1.2cm}\times
\frac{[u-v_{l-1}+(P+h)_{j,l}+\frac{l-N}{2}-1]
[1]}{[u-v_{l-1}+\frac{l-N}{2}]
[P_{j,l}+h_{j,l}-1]}
\prod_{m=j}^{l-2}
\frac{[v_{m+1}-v_{m}+(P+h)_{j,m+1}-\frac{1}{2}][1]}
{[v_{m+1}-v_{m}+\frac{1}{2}][P_{j,m+1}+h_{j,m+1}
]},\lb{fjl}\\
E_{l,j}^+(z)&=&a_{j,l}^*\oint_{\T^{l-j}
}
\prod_{m=j}^{l-1}\frac{dt_m}{2\pi i t_m}
E_{l-1}(t_{l-1})E_{l-2}(t_{l-2})\cdots E_{j}(t_{j})
\nonumber\\
&&\hspace{-1cm}\times
\frac{[u-v_{l-1}-P_{j,l}+\frac{l-N}{2}-
\frac{c}{2}+1]^*[1]^*}
{[u-v_{l-1}+\frac{l-N}{2}-\frac{c}{2}]^*[P_{j,l}-1]^*
}
\prod_{m=j}^{l-2}
\frac{[v_{m+1}-v_{m}-P_{j,m+1}+\frac{1}{2}]^*[1]^*}
{[v_{m+1}-v_{m}-\frac{1}{2}]^*[P_{j,m+1}-1]^*},\lb{elj}
\end{eqnarray}  
where $z=q^{2u}$, $t_a=q^{2v_a}$ $(a=j,j+1,\cdots,l-1)$ and 
\be
&&\T^{l-j}=\{t\in \C^{l-j}\ |\ |t_j|=\cdots=|t_{l-1}|=1 \}. 
\en
The constants $a_{j,l}$ and $a_{j,l}^*$ are chosen to 
satisfy
\begin{eqnarray}
&&-\frac{\varrho
\ a_{j,l}a_{j,l}^*}{q-q^{-1}}\frac{ [1]}{[0]'}=1.
\end{eqnarray}
We call $F_{j,l}^+(z), E_{l,j}^+(z), (1\leq j<l \leq N)$
and $K_j^+(z)\ (j=1,\cdots,N)$ the half currents.
\end{dfn}

\begin{dfn}\lb{def:Lop} By using the half currents, we define the $L$-operator
$\widehat{L}^+(z) \in {\rm End}({\mathbb{C}}^N)
\otimes U_{q,p}(\slnh)$
as follows.
\begin{eqnarray}
&&\widehat{L}^+(z)=
\left(\begin{array}{ccccc}
1&F_{1,2}^+(z)&F_{1,3}^+(z)&\cdots&F_{1,N}^+(z)\\
0&1&F_{2,3}^+(z)&\cdots&F_{2,N}^+(z)\\
\vdots&\ddots&\ddots&\ddots&\vdots\\
\vdots&&\ddots&1&F_{N-1,N}^+(z)\\
0&\cdots&\cdots&0&1
\end{array}\right)\left(
\begin{array}{cccc}
K^+_1(z)&0&\cdots&0\\
0&K^+_2(z)&&\vdots\\
\vdots&&\ddots&0\\
0&\cdots&0&K^+_{N}(z)
\end{array}
\right)
\nn\\
&&\qquad\qquad\qquad\qquad\qquad\qquad\qquad\times
\left(
\begin{array}{ccccc}
1&0&\cdots&\cdots&0\\
E^+_{2,1}(z)&1&\ddots&&\vdots\\
E^+_{3,1}(z)&
E^+_{3,2}(z)&\ddots&\ddots&\vdots\\
\vdots&\vdots&\ddots&1&0\\
E^+_{N,1}(z)&E^+_{N,2}(z)
&\cdots&E^+_{N,N-1}(z)&1
\end{array}
\right).\lb{def:lhat}
\end{eqnarray}
\end{dfn}
We conjecture that the $L$-operator satisfies the following dynamical $RLL$-relation\cite{JKOStg}.
\noindent
\begin{conj}\lb{main}\cite{KK03}
\begin{eqnarray}
R^{+(12)}(z_1/z_2,\Pi)\widehat{L}^{+(1)}(z_1)
\widehat{L}^{+(2)}(z_2)=
\widehat{L}^{+(2)}(z_2)
\widehat{L}^{+(1)}(z_1)
R^{+*(12)}(z_1/z_2,\Pi^*).\label{thm:RLL}
\end{eqnarray}
\end{conj}
\noindent

\subsection{The $H$-Hopf algebroid $U_{q,p}(\slnh)$}\lb{Hopfalgebroid}

Let $\cA$ be a complex associative algebra, $\cH$ be a finite dimensional commutative subalgebra of 
$\cA$, and $\cM_{\cH^*}$ be the 
field of meromorphic functions on $\cH^*$ the dual space of $\cH$. 

\begin{dfn}[$\cH$-algebra\cite{EV}] 
An $\cH$-algebra is an associative algebra $\cA$ with 1, which is bigraded over 
$\cH^*$, $\ds{\cA=\bigoplus_{\alpha,\beta\in \cH^*} \cA_{\al,\beta}}$, and equipped with two 
algebra embeddings $\mu_l, \mu_r : \cM_{\cH^*}\to \cA_{0,0}$ (the left and right moment maps), such that 
\be
\mu_l(\hf)a=a \mu_l(T_\al \hf), \quad \mu_r(\hf)a=a \mu_r(T_\beta \hf), \qquad 
a\in \cA_{\al,\beta},\ \hf\in \cM_{\cH^*},
\en
where $T_\al$ denotes the automorphism $(T_\al \hf)(\la)=\hf(\la+\al)$ of $\cM_{\cH^*}$.
\end{dfn}

Let $\cA$ and $\cB$ be two $\cH$-algebras. The tensor product $\cA {\widetilde{\otimes}}\cB$ is the $\cH^*$-bigraded vector space with 
\be
 (\cA {\widetilde{\otimes}}\cB)_{\al\beta}=\bigoplus_{\gamma\in\H^*} (\cA_{\al\gamma}\otimes_{\cM_{\cH^*}}\cB_{\gamma\beta}),
\en
where $\otimes_{\cM_{\cH^*}}$ denotes the usual tensor product 
modulo the following 
relation.
\bea
\mu_r^\cA(\hf) a\otimes b=a\otimes\mu_l^\cB(\hf) b, \qquad a\in \cA, 
b\in \cB, \hf\in \cM_{\cH^*}.\lb{AtotB}
\ena
The tensor product $\cA {\widetilde{\otimes}}\cB$ is again an $\cH$-algebra with the multiplication $(a\otimes b)(c\otimes d)=ac\otimes bd$ and the moment maps 
\be
\mu_l^{\cA {\widetilde{\otimes}}\cB} =\mu_l^\cA\otimes 1,\qquad \mu_r^{\cA {\widetilde{\otimes}}\cB} =1\otimes \mu_r^\cB.
\en

\begin{prop}\cite{Konno09,Konno16}
$U=U_{q,p}(\slnh)$  is an $H$-algebra by 
\be
&&U=\bigoplus_{\al,\beta\in H^*}U_{\al,\beta}\\
&&U_{\al,\beta}=\left\{x\in U \left|\ q^{P+h}x q^{-(P+h)}=q^{<\al,P+h>}x,\quad q^{P}x q^{-P}=q^{<\beta,P>}x\ \forall P+h, P\in H\right.\right\}
\en
and $\mu_l, \mu_r : \FF \to U_{0,0}$ defined by 
\be
&&\mu_l(\hf)=f(P+h,p)\in \FF[[p]],\qquad \mu_r(\hf)=f(P,p^*)\in \FF[[p]].
\en
\end{prop}
We regard ${T}_{\al}=e^\al\in \C[\cR_Q]$ as 
the shift operator $\cM_{{H}^*}\to \cM_{{H}^*}$ 
\be
({T}_{\al}\widehat{f})=e^{\al}f(P,p^*)e^{-\al}={f}(P+<\al,P>,p^*).
\en
Hereafter we abbreviate 
$f(P+h,p)$ and $f(P,p^*)$ as $f(P+h)$ and
 $f^*(P)$, respectively.


We also consider the $H$-algebra of the shift operators 
\be
&&\cD=\{\ \sum_\al \widehat{f}_\al{T}_{\al}\ |\ \widehat{f}_\al\in {M}_{H^*}, 
\al\in \cR_Q\ \},\\
&&\cD_{\ha,\ha}=\{\ \widehat{f}{T}_{-\ha}\ \},\quad \cD_{\ha,{\beta}}=0\ 
(\al\not=\beta),\\
&&\mu_l^{\cD}(\widehat{f})=
\mu_r^{\cD}(\widehat{f})=\widehat{f}{T}_0 \qquad \widehat{f}\in {M}_{H^*}.
\en
Then we have the $H$-algebra isomorphism 
\bea
 U\cong U\tot\cD\cong \cD\tot U. \lb{Diso}
\ena

We define two $H$-algebra homomorphisms, the co-unit $\vep : U\to \cD$ and the co-multiplication $\Delta : U\to U \widetilde{\otimes}U$ 
as well as the algebra antihomomorphism $S:U\to U$ by
\bea
&&\vep(\hL^+_{i,j}(z))=\delta_{i,j}{T}_{Q_{\ep_i} }\quad (n\in \Z),
\qquad \vep(e^Q)=e^Q,\lb{counitUqp}\\
&&\vep(\mu_l({\hf}))= \vep(\mu_r(\hf))=\widehat{f}T_0, \lb{counitf}\\
&&\Delta(\hL^+_{i,j}(z))=\sum_{k}\hL^+_{i,k}(z)\widetilde{\otimes}
\hL^+_{k,j}(z),\lb{coproUqp}\\
&&\Delta(e^{Q})=e^{Q}\tot e^{Q},\\
&&\Delta(\mu_l(\hf))=\mu_l(\hf)\widetilde{\otimes} 1,\quad \Delta(\mu_r(\hf))=1\widetilde{\otimes} \mu_r(\hf),\lb{coprof}\\
&&S(\hL^+_{ij}(z))=(\hL^+(z)^{-1})_{ij},\\
&&S(e^{Q})=e^{-Q}, \quad S(\mu_r(\hat{f}))=\mu_l(\hat{f}),\quad S(\mu_l(\hat{f}))=\mu_r(\hat{f}).
\ena
Then the set  $(U_{q,p}(\slnh),H,{\cM}_{H^*},\mu_l,\mu_r,\Delta,\vep,S)$ becomes an $H$-Hopf algebroid\cite{EV,KR,Konno09,Konno16}.

\section{Representations}\lb{App:rep}
\subsection{Dynamical representations}
Let us consider a vector space $\hV$ over $\FF=\cM_{H^*}$, which is  
${H}$-diagonalizable, i.e.  
\be
&&\hV=\bigoplus_{\la,\nu\in {H}^*}\hV_{\la,\nu},\ \hV_{\la,\nu}=\{ v\in \hV\ |\ q^{P+h}\cdot v=q^{<\la,P+h>} v,\ q^{P}\cdot v=q^{<\nu,P>} v\ \forall 
P+h, P\in 
{H}\}.
\en
Let us define the $H$-algebra $\cD_{H,\hV}$ of the $\C$-linear operators on $\hV$ by
\be
&&\cD_{H,\hV}=\bigoplus_{\al,\beta\in {H}^*}(\cD_{H,\hV})_{\al,\beta},\\
&&\hspace*{-10mm}(\cD_{H,\hV})_{\al,\beta}=
\left\{\ X\in \End_{\C}\hV\ \left|\ 
\mmatrix{ f(P+h)X=X f(P+h+<\alpha,P+h>),\cr 
f(P)X=X f(P+<\beta,P>)\cr
 f(P), f(P+h)\in \FF,\ X\cdot\hV_{\la,\mu}\subseteq 
 \hV_{\la+\al,\mu+\beta}\cr}  
 \right.\right\},\\
&&\mu_l^{\cD_{H,\hV}}(\widehat{f})v=f(<\la,P+h>,p)v,\quad 
\mu_r^{\cD_{H,\hV}}(\widehat{f})v=f(<\nu,P>,p^*)v,\quad \widehat{f}\in \FF,
\ v\in \hV_{\la,\nu}.
\en
\begin{dfn}\cite{EV,KR,Konno09} 
We define a dynamical representation of $U_{q,p}(\slnh)$ on $\hV$ to be  
 an $H$-algebra homomorphism ${\pi}: U_{q,p}(\slnh) 
 \to \cD_{H,\hV}$. By the action of $U_{q,p}(\slnh)$ we regard $\hV$ as a 
$U_{q,p}(\slnh)$-module. 
\end{dfn}

\begin{dfn}
For $k\in \C$, we say that a $U_{q,p}(\slnh)$-module has  level $k$ if $c$ act 
as the scalar $k$ on it.  
\end{dfn}


\begin{dfn}
Let ${\cal H}$, ${\cal N}_+, {\cal N}_-$ be the subalgebras of 
$U_{q,p}(\slnh)$ 
generated by 
$c, d, 
K^\pm_{i}\ (i\in I)$, by $\al_{i,n}\ (i\in I, n\in \Z_{>0})$,  $e_{i, n}\ (i\in I, n\in \Z_{\geq 0})$  
$f_{i, n}\ (i\in I, n\in \Z_{>0})$ and by $\al_{i,-n}\ (i\in I, n\in \Z_{>0}),\ e_{i, -n}\ (i\in I, n\in \Z_{> 0}),\ f_{i, -n}\ (i\in I, n\in \Z_{\geq 0})$, respectively.   
\end{dfn}

\begin{dfn}\lb{def:levelkrep}
For $k\in\C$, $\la\in \h^*$ and $\nu\in H^*$, 
a (dynamical) $U_{q,p}(\slnh)$-module $\hV(\la,\nu)$ is called the 
level-$k$ highest weight module with the highest weight $(\la,\nu)$, if there exists a vector 
$v\in \hV(\la,\nu)$ such that
\be
&&\hV(\la,\nu)=U_{q,p}(\slnh)\cdot v,\qquad \cN_+\cdot v=0,\\
&&c\cdot v=kv, 
\quad  f({P})\cdot v =f({<\nu,P>})v,\quad f({P+h})\cdot v =f({<\la,P+h>})v.
\en
\end{dfn}

\subsection{The $N$-dimensional dynamical evaluation representation}\lb{vecrep}


Let  $\displaystyle{\hV=\bigoplus_{\mu=1}^{ N} \FF v_\mu\otimes 1}$ and set $\hV_z=\hV[z,z^{-1}]$. Let $e^{Q_{\al}}\in \C[\cR_Q]$ act on $f(P_\beta)v\otimes 1$  by $e^{Q_{\al}}(f(P_\beta)v\otimes 1)=f(P_\beta-(\al,\beta))v\otimes 1$. 
\begin{thm}\lb{vecRep}
Let $E_{j,k}\ (1\leq j,k\leq N)$ denote the matrix units such that 
$E_{j,k}v_\mu=\delta_{k,\mu}v_j$. 
The following gives the $N$-dimensional dynamical evaluation representation of $U_{q,p}(\slnh)$  
on $\hV_z$.  
\be
&&\pi_z(c)=0,\quad \pi_z(d)=-z\frac{d}{dz},\\
&&\pi_z(\al_{j,m})=\frac{[m]_q}{m}(q^{j-N+1}z)^n(q^{-m}E_{j,j}-q^mE_{j+1,j+1}),\\
 && \pi_z(e_j(w))=\frac{(pq^2;p)_{\infty}}{(p;p)_{\infty}}E_{j,j+1}
     \delta\left(q^{j-N+1}{z}/{w}\right) 
    e^{-Q_{\al_j}}, \\
 && \pi_z(f_j(w))=\frac{(pq^{-2};p)_{\infty}}{(p;p)_{\infty}}E_{j+1,j}
     \delta\left(q^{j-N+1}{z}/{w}\right), \\
 && \pi_z(\psi_j^{+}(w,p))=q^{-\pi(h_j)}e^{-Q_{\al_j}}
     \frac{\Theta_p(q^{-j+N-1+2\pi(h_j)} \frac{w}{z})}{\Theta_p(q^{-j+N-1} {w}/{z})}
    \qquad (1\leq j\leq N-1),\\
 && \pi_z(\psi_j^{-}(w,p))=q^{\pi(h_j)}e^{-Q_{\al_j}}
     \frac{\Theta_p(q^{j-N+1-2\pi(h_j)} \frac{z}{w})}{\Theta_p(q^{j-N+1} {z}/{w})}. 
\en
Here $\pi(h_j)=E_{j,j}-E_{j+1,j+1}$.
\end{thm}

Combining the formulas in Definition \ref{modellcurrents}, \ref{def:halfcurrents}, \ref{def:Lop} and Theorem \ref{vecRep}, we obtain 
\begin{cor}\lb{repLR}
\be
&&\pi_z({\hL^+_{i,j}(w)})_{k,l}=R^+(w/z,\Pi^*)_{ik}^{jl}.
\en
\end{cor}

\subsection{The level-1 representation}\lb{levelonerep}
Next we consider level-1 ($c=1$) representation of $U_{q,p}(\slnh)$.
We mainly follow the work \cite{FKO}.  

It is convenient to  extend the  root lattice ${\cQ}$ by adding the elements $\zeta_j\ (j=1,\cdots,N-1)$ in Definition \ref{extragen}. 
Let us set $\widehat{\alpha}_j=\alpha_j+\zeta_j$ and consider $\widehat{\cQ}=\oplus_j\Z \widehat{\alpha}_j$.   
We define the extended group algebra 
 $\C[\widehat{\cQ}]$ with assuming the following central extension. 
\be
&&e^{ \alpha_i}e^{ \alpha_j}=(-1)^{(\alpha_i, \alpha_j)} e^{\alpha_j}e^{ \alpha_i}. 
\en

For  $\omega=\sum_jc_j\bep_j\in \h^*$, we also set 
$\zeta_{\omega}=\sum_jc_j \zeta_{\bep_j}$ and $\widehat{\omega}=\omega+\zeta_{\omega}$. 
Set also $\widehat{\Lambda}_0=\Lambda_0$.

Let $\Lambda_0$ and $\Lambda_a=\Lambda_0+\bar{\Lambda}_a \ (a=1,\cdots,N-1)$ be the fundamental weights of $\slnh$. For generic 
$\nu\in \h^*$, we set 
\be
\hV(\Lambda_{a}+\nu,\nu)&=&\FF\otimes_\C (\F_{\al,1}\otimes e^{\widehat{\Lambda}_a}\C[\widehat{\cQ}])\otimes 
e^{Q_{\bar{\nu}}}\C[{\cR}_Q],
\en 
where $\F_{\al,1}=\C[\{\alpha_{j,-m}\ (j=1,\cdots,N-1,\ m\in \N_{>0})\}]$. 
Then we have the following decomposition.
\be
\hV(\Lambda_{a}+\nu,\nu)&=&\bigoplus_{\xi, \eta\in \cQ}
\F_{a,\nu}(\xi,\eta),
\en
where 
\bea
&&\F_{a,\nu}(\xi,\eta)=\FF\otimes_\C(\F_{\al,1}\otimes e^{\widehat{\Lambda}_a+\widehat{\xi}})\otimes  e^{Q_{{\nu}+\eta}}.\lb{repsp}
\ena

\begin{thm}\lb{levelone}\cite{FKO} 
The spaces $\hV(\Lambda_{a}+\nu,\nu)$ ($a=0,\cdots,N$) give the level-1 irreducible $U_{q,p}(\slnh)$-modules with the highest weight $(\Lambda_a+\nu,\nu)$, where 
the highest weight vectors are given by $1\otimes  e^{\Lambda_a}\otimes e^{Q_{\bar{\nu}}}$. 
The  action of the elliptic currents is given by 
\begin{eqnarray}
 && E_j(z) =\; :\exp \left\{ -\sum_{n \ne 0} \frac{1}{[n]_q} \alpha_{j,n}z^{-n} \right\}:
 e^{\widehat{\alpha}_j}e^{-Q_{\alpha_j}}z^{h_{\alpha_j}+1}(q^{N-j}z)^{-\frac{P_{\alpha_j}-1}{r^*}},  \\
 && F_j(z) =\; :\exp \left\{ \sum_{n \ne 0}\frac{1}{[n]_q} \alpha'_{j,n}z^{-n}\right\}: e^{-\widehat{\alpha}_j}z^{-h_{\alpha_j}+1}
(q^{N-j}z)^{\frac{(P+h)_{\alpha_j}-1}{r}}, \lb{level1Fj}\\
&&\hspace{7cm} (1\leq j\leq N-1)\nn
\end{eqnarray}
together with $H_j^{\pm}(z), \,K_{j}^{+}(z)$  in Sec.\ref{MEC} and 
\be
&&\hd=d+\frac{1}{2r^*}\sum_{j=1}^N(P_j+2)P^j-\frac{1}{2r}\sum_{j=1}^N((P+h)_j+2)(P+h)^j,\\
&&d=-\frac{1}{2}\sum_{j=1}^{N-1}h_jh^j-\sum_{j=1}^{N-1}\sum_{m\in \Z_{>0}}\frac{m^2}{[m]}\frac{1-p^{*m}}{1-p^m}q^m\al_{j,-m}A^j_m.
\en
In \eqref{level1Fj} we set
 $\ds{\alpha'_{j,n}=\frac{1-p^{*n}}{1-p^n}q^{n} \alpha_{j,n}}$. 
\end{thm}

\section{Proof of Proposition \ref{shuffleMp}}\lb{shufflealg}
Let $\la, \la'\in \N^N$, $|\la|=m$, $|\la'|=n$ and consider 
 $I=I_{\mu_1\cdots\mu_m}=(I_1,\cdots,I_N)\in \cI_{\la}$ and $I'=I_{\mu'_1\cdots\mu'_m}=(I'_1,\cdots,I'_N)\in \cI_{\la'}$. 
 For each $I'_l=\{i'_{l,1},\cdots,i'_{l,\la'_l}\}$ $(l=1,\cdots,N)$, let us set $\widetilde{I'}_l=\{m+i'_{l,1},\cdots,m+i'_{l,\la'_l}\}$. 
Then define $I+I'=((I+I')_1,\cdots,(I+I')_N)\in \cI_{\la+\la'}$ by $(I+I')_l=I_l\cup \widetilde{I'}_l\ (l=1\cdots,N)$.  
 
Let us consider the $m$- and $n$-point functions $\phi_{\mu_1\cdots\mu_m}(z_1,\cdots,z_m)$ and  
$\phi_{\mu'_1,\cdots,\mu'_{n}}(z'_1,\cdots,z'_{n})$. 
Their composition gives the $m+n$-point function 
\bea
&&\hspace{-0.5cm}\phi_{\mu_1,\cdots,\mu_{m+n}}(z_1,\cdots,z_{m+n})=\phi_{\mu_1\cdots\mu_m}(z_1,\cdots,z_m)\phi_{\mu'_1,\cdots,\mu'_{n}}(z'_1,\cdots,z'_{n}), \lb{phiphimphin}
\ena
where we set $z_{m+k}:=z'_{k}$ and $ \mu_{m+k}:=\mu'_k$ $(k=1,\cdots,n)$.

On the otherhand, from Theorem \ref{prodPhi} we have
\bea
&&\phi_{\mu_1,\cdots,\mu_{m}}(z_1,\cdots,z_{m})=\oint_{\T^{M}} \underline{dt}\ \widetilde{\Phi}(t,z)
\omega_{\mu_1\cdots \mu_{m}}(t,z,\Pi_I),
\lb{phim}\\
&&\phi_{\mu'_1,\cdots,\mu'_{n}}(z'_1,\cdots,z'_{n})=\oint_{\T^{M'}} \underline{dt'}\ \widetilde{\Phi}(t',z')
\omega_{\mu'_1\cdots \mu'_{n}}(t',z',\Pi'_{I'}),\lb{phin}\\
&&\phi_{\mu_1,\cdots,\mu_{m+n}}(z_1,\cdots,z_{m+n})=\oint_{\T^{{M}+M'}} \underline{d\tilde{t}}\ \widetilde{\Phi}(\tilde{t}, z\cup z')\omega_{\mu_1\cdots \mu_{m+n}}(\tilde{t},z\cup z',\widetilde{\Pi}).\lb{phimpn}
\ena
Here we set $z=(z_1,\cdots, z_m)$,  $z'=(z'_{1},\cdots,z'_{n})$, $M=\sum_{j=1}^{N-1}(N-j)\la_j$, $M'=\sum_{j=1}^{N-1}(N-j)\la'_j$, 
 $t=(t^{(1)}_1,\cdots,t^{(1)}_{\la^{(1)}},\cdots,t^{(N-1)}_1,\cdots, t^{(N-1)}_{\la^{(N-1)}})$, 
$t'=(t^{'(1)}_1,\cdots,t^{'(1)}_{\la^{'(1)}},\cdots,$ $ t^{'(N-1)}_1,\cdots, t^{'(N-1)}_{\la^{'(N-1)}})$ 
and $\tilde{t}=(\tilde{t}^{(1)}_1,\cdots,\tilde{t}^{(1)}_{{\la}^{(1)}+{\la'}^{(1)}},\cdots,\tilde{t}^{(N-1)}_1,\cdots, \tilde{t}^{(N-1)}_{{\la}^{(N-1)}
+{\la'}^{(N-1)}})$ and $\widetilde{\Pi}=\{ \widetilde{\Pi}_{\mu_{k},j}$ $(k=1,\cdots,m+n,\ j=\mu_k+1,\cdots,N)\}$.  

Substituting \eqref{phim}-\eqref{phimpn} into \eqref{phiphimphin}, one can obtain a relation among 
 the weight functions $\omega_{\mu_1\cdots \mu_{m}}(t,z,\Pi_I)$, $\omega_{\mu'_1\cdots \mu'_{n}}(t',z',\Pi'_{I'})$ and 
 $\omega_{\mu_1\cdots \mu_{m+n}}(\tilde{t},z\cup z',\widetilde{\Pi})$. 
 
\begin{dfn}
For the weight functions  $\omega_{\mu_1\cdots\mu_m}(t,z,\Pi)$ and $\omega_{\mu'_1\cdots\mu'_n}(t',z',\Pi')$, 
we define the  $\star$-product as follows.
\be
&&(\omega_{\mu_1\cdots\mu_m}\star \omega_{\mu'_1\cdots\mu'_n})(t\cup t' 
,z\cup z',\Pi_{I+I'})\\
&&\hspace{-1cm}:=\frac{1}{\prod_{l=1}^{N-1}\la^{(l)}!\la^{'(l)}!}{\rm Sym}_{t^{(1)}}\cdots {\rm Sym}_{t^{(N-1)}} \left[
\omega_{\mu_1\cdots\mu_m}(t,z,\Pi_I q^{-2\sum_{j=1}^n<\bep_{\mu_j'},h>})\ \omega_{\mu'_1\cdots\mu'_n}(t',z',\Pi'_{I'})\ \widetilde{\Xi}(t,t',z,z')
\right], 
\en
where   ${\Pi}_{I+I'}=\{\Pi_{\mu_{k},j}$ $(k=1,\cdots,m+n,\ j=\mu_k+1,\cdots,N)$ 
with ${\Pi}_{\mu_{m+k},j}=\Pi'_{\mu'_k,j}$ $(k=1,\cdots,n,\ j=\mu'_k+1,\cdots,N)$, and   
\be
&&\widetilde{\Xi}(t,t',z,z')=\mu_{m,n}^+(z,z')\prod_{l=1}^{N-1}
\prod_{a=1}^{\la^{(l)}}\left(\prod_{b=1}^{\la^{'(l+1)}}\frac{[{v'_b}^{(l+1)}-v^{(l)}_a-\frac{1}{2}]}{[{v'_b}^{(l+1)}-v^{(l)}_a+\frac{1}{2}]}
\prod_{c=1}^{\la^{'(l)}}\frac{[{v_c'}^{(l)}-v^{(l)}_a+1]}{[{v'_c}^{(l)}-v^{(l)}_a]}\right),\\
&&\mu_{m,n}^+(z,z')=\prod_{k=1}^m\prod_{l=1}^n(-)^{\frac{N-1}{N}}z_k^{\frac{r-1}{r}\frac{N-1}{N}}\frac{\Gamma(q^2z_l'/z_k;p,q^{2N})}{\Gamma(q^{2N}z_l'/z_k;p,q^{2N})}.
\en
\end{dfn}

\begin{prop}\lb{shuffleprod}
Under the identification 
\bea
&&\tilde{t}^{(l)}_a={t}^{(l)}_a\quad (a=1,\cdots,\la^{(l)}),\qquad  \tilde{t}^{(l)}_{\la^{(l)}+b}={t'}^{(l)}_b \quad (b=1,\cdots,{\la'}^{(l)})\lb{identifyt}
\ena
 for $l=1,\cdots,N-1$ and $\Pi_{I+I'}=\widetilde{\Pi}$, we obtain  
\be
&&(\omega_{\mu_1\cdots\mu_m}\star \omega_{\mu'_1\cdots\mu'_n})(t\cup t' 
,z\cup z',\Pi_{I+I'})=\omega_{\mu_1\cdots \mu_{m+n}}(\tilde{t},z\cup z',\widetilde{\Pi}).
\en
\end{prop}

\noindent
{\it Proof.} \ 
%
Substituting \eqref{phim} and \eqref{phin} into the right hand side of \eqref{phiphimphin}, we obtain
\be
&&\phi_{\mu_1\cdots\mu_{m+n}}(z_1,\cdots,z_{m+n})\nn\\
&&=\oint_{\T^{M}} \underline{dt}\oint_{\T^{M'}} \underline{dt'}\ \widetilde{\Phi}(t,z)\omega_{\mu_1\cdots\mu_m}(t,z,\Pi_I)
\widetilde{\Phi}(t',z')\omega_{\mu'_1\cdots\mu'_m}(t',z',\Pi'_{I'})\nn\\
&&=\oint_{\T^{M}} \underline{dt}\oint_{\T^{M'}} \underline{dt'}\ \widetilde{\Phi}(t,z)\widetilde{\Phi}(t',z')\omega_{\mu_1\cdots\mu_m}(t,z,\Pi_I q^{-2\sum_{j=1}^n<\bep_{\mu'_j},h>})\omega_{\mu'_1\cdots\mu'_m}(t',z',\Pi'_{I'}).
\en
Furthermore from  \eqref{u8} and \eqref{PhiNFNm1} we obtain
\be
&& \widetilde{\Phi}(t,z)\widetilde{\Phi}(t',z')= {\Upsilon}(t, t',z, z')\widetilde{\Xi}(t,t',z,z')
\en
where
\be
&&{\Upsilon}(t, t',z, z')\\
&&=:\Phi_N(z_1)\cdots \Phi_N(z_{m+n}): 
:F_{N-1}(t^{(N-1)}_{1})\cdots F_{N-1}(t^{(N-1)}_{\la^{(N-1)}})\cdots F_{N-1}(t^{'(N-1)}_{1})\cdots F_{N-1}(t^{'(N-1)}_{\la^{'(N-1)}}):\\
&&\qquad\qquad\qquad  \cdots :F_1(t_1^{(1)})\cdots F_{1}(t^{(1)}_{\la^{(1)}})\cdots  F_1(t_1^{'(1)})\cdots F_{1}(t^{'(1)}_{\la^{'(1)}}): 
\nn\\
&&\times \prod_{1\leq k<l\leq m}<\Phi_N(z_k)\Phi_N(z_l)>^{Sym}\prod_{1\leq k<l\leq n}<\Phi_N(z'_k)\Phi_N(z'_l)>^{Sym}
\prod_{k=1}^{m}\prod_{l=1}^{n}<\Phi_N(z_k)\Phi_N(z'_l)>^{Sym}\\
&&\hspace{-0.5cm}\times\prod_{l=1}^{N-1}\left(\prod_{1\leq a<b\leq \la^{(l)}}<F_l(t^{(l)}_a)F_l(t^{(l)}_b)>^{Sym}\prod_{1\leq a<b\leq \la^{'(l)}}<F_l(t^{'(l)}_a)F_l(t^{'(l)}_b)>^{Sym}\prod_{a=1}^{\la^{(l)}}\prod_{b=1}^{\la^{'(l)}}<F_l(t^{(l)}_a)F_l(t^{(l)}_b)>^{Sym}
\right).
\en
Then it turns out that  $ {\Upsilon}(t, t',z, z')$ coinsides with  $\widetilde{\Phi}(\tilde{t}, z\cup z')$ under the identification \eqref{identifyt}. 
Hence we have
\be
&&\phi_{\mu_1\cdots\mu_{m+n}}(z_1,\cdots,z_{m+n})\nn\\
&&=\oint_{\T^{{M}+M'}} \underline{d\tilde{t}}\ \widetilde{\Phi}(\tilde{t},z\cup z')\omega_{\mu_1\cdots\mu_m}(t,z,\Pi_I q^{-2\sum_{j=1}^n<\bep_{\mu'_j},h>})\omega_{\mu'_1\cdots\mu'_m}(t',z',\Pi'_{I'})\widetilde{\Xi}(t,t',z,z')\\
&&=\oint_{\T^{{M}+M'}} \underline{d\tilde{t}}\ \widetilde{\Phi}(\tilde{t},z\cup z')
\frac{1}{\prod_{l=1}^{N-1}\la^{(l)}!\la^{'(l)}!}{\rm Sym}_{t^{(1)}}\cdots {\rm Sym}_{t^{(N-1)}} \left[
\omega_{\mu_1\cdots\mu_m}(t,z,\Pi_I q^{-2\sum_{j=1}^n<\bep_{\mu'_j},h>})\right.\\
&&\left.\qquad\qquad\qquad\qquad\qquad\qquad\qquad \times\omega_{\mu'_1\cdots\mu'_m}(t',z',\Pi'_{I'})\widetilde{\Xi}(t,t',z,z')
\right].
\en
The last equality follows from the symmetry of $\widetilde{\Phi}(\tilde{t}, z\cup z')$ under the action of 
$\sigma\in \gS_{\widetilde{\la}^{(1)}}\times \cdots \times \gS_{\widetilde{\la}^{(N-1)}}$ on $\tilde{t}$. 
Comparing this with \eqref{phimpn} we obtain the desired result. 
\qed

Proposition \ref{shuffleMp} is a direct consequence of this proposition. 

\end{appendix}

\renewcommand{\baselinestretch}{0.7}


\begin{thebibliography}{}
\bibliographystyle{unsrt}


\bibitem{AO}
M. Aganagic and A. Okounkov, 
\newblock Elliptic Stable Envelopes, 
Preprint 2016, arXiv:1604.00423.

\bibitem{AFO}
M. Aganagic E.Frenkel and A. Okounkov, 
\newblock Quantum $q$-Langlands Correspondence, 
Preprint 2017, arXiv:1701.03146.

\bibitem{AKOS}
H. Awata, H. Kubo, S. Odake, and J. Shiraishi, 
\newblock Quantum $W_N$ algebras and Macdonald
Polynomials,
{\it  Comm. Math. Phys.} {\bf 179} (1996) 401--416.

\bibitem{Drinfeld}
 V.G. Drinfeld, A New Realization of Yangians and Quantized
Affine Algebras. {\it Soviet Math. Dokl.}{\bf  36} (1988) 212-216.


\bibitem{EV}
{P. Etingof and A. Varchenko},
{{Solutions of the Quantum Dynamical Yang-Baxter Equation and Dynamical Quantum Groups}},
{\it Comm.Math.Phys.},
{\bf 196},
 1998, 591--640; 
  {{Exchange Dynamical Quantum Groups}},
  {\it Comm.Math.Phys.},
  {\bf 205},
   (1999) 19--52.


\bibitem{FKO}
     {R. M. Farghly, H. Konno and K. Oshima},
     {{Elliptic Algebra $U_{q,p}(\gh)$ and Qunatum $Z$-algebras}},
     Algebras and Representation Theory (2014) June, DOI 10.1007/s10468-014-9483-x. 


\bibitem{FF}
B. Feigin and E. Frenkel, 
\newblock Quantum $W$-algebras and Elliptic Algebras,
{\it Comm. Math.Phys.} {\bf 178} (1996) 653--678. 


\bibitem{FJMOP}
B.Feigin, M.Jimbo, T.Miwa, A.Odesskii and Y.Pugai,
\newblock Algebra of Screening Operators for the Deformed $W_n$ Algebra, 
{\it Comm.Math.Phys.} {\bf 191} (1998) 501--541.  


\bibitem{FT}
B.Feigin and A.Tsymbaliuk,
\newblock Bethe Subalgebras of $U_q(\widehat{\gl}_n)$ via Shuffle Algebras,
{\it Selecta Math. (N.S.)} {\bf 22} (2016), no. 2, 979--1011.





\bibitem{Felder1}
G.Felder, 
\newblock Elliptic Quantum Groups, in Proceedings of the XIth ICMP, Paris 1994, 
Intern. Press, Cambridge, MA, (1995) 211--218.

\bibitem{Felder2}
G.Felder, 
\newblock Conformal Field Theory and Integrable Systems Associated to Elliptic Curves. 
Proceedings of the International Congress of Mathematicians, Vol. 1, 2 (Z\"urich, 1994), 1247--1255, 
Birkh\" auser, Basel, 1995.

 
\bibitem{FTV1}
G.Felder, V.Tarasov  and A.Varchenko, 
\newblock  Solutions of the Elliptic QKZB Equations and Bethe ansatz I, 
in Topics in Singularity Theory, V.I.Arnold’s 60th Anniversary Collection, Advances in the Mathematical Sciences, 
AMS Translations, Series 2, v.180,(1997) 45--76.
  
\bibitem{FTV2}
G.Felder, V.Tarasov  and A.Varchenko, 
\newblock Monodromy of Solutions of the Elliptic Quantum Knizhnik-Zamolodchikov-Bernard Difference Equations, 
{\it Int. J. Math. } {\bf10} (1999) 943--975.



\bibitem{FRV}
G.Felder, R.Rim\' anyi  and A.Varchenko, 
\newblock Elliptic Dynamical Quantum Groups and Equivariant Elliptic Cohomology, 
Preprint (2017), arXiv:1702.08060



\bibitem{FJMMN}
O.Foda, M.Jimbo, T.Miwa, K.Miki and A.Nakayashiki,
\newblock Vertex Operators in Solvable Lattice Models,
{\it J.Math.Phys.} {\bf 35} (1994) 13--46.





\bibitem{GRTV}
V.Gorbounov, R.Rim\' anyi, V.Tarasov  and A.Varchenko, 
\newblock Cohomology of the Cotangent Bundle of a Flag Variety as a Yangian Bethe Algebra,
 {\it J.Geom.Phys.},
  {\bf 74} (2013) 56--86.



\bibitem{JMO}
 	 {M. Jimbo and T. Miwa  and M. Okado},
 	 {{Solvable Lattice Models Related to the Vector Representation of Classical Simple Lie Algebras}},
  	 {{\it Comm. Math. Phys.}}
   {\bf 116}
  {(1988)}
  	 {507--525}.


\bibitem{JM}
M. Jimbo and T. Miwa, 
 Algebraic Analysis of Solvable Lattice Models.
 {\it Conference Board of the Math. Sci., Regional Conference Series
 in Mathematics}{\bf  85} (1995) and references therein.
 
 
\bibitem{JKOStg}
 	 {M. Jimbo, H. Konno, S. Odake and J. Shiraishi},
  	 {{Quasi-Hopf Twistors for Elliptic Quantum Groups}},
  	 {{\it Transformation Groups}}
  	 {\bf 4} {(1999)} {303--327}.


\bibitem{JKOS}
 	 {M. Jimbo, H. Konno, S. Odake and J. Shiraishi},
  	 {{Elliptic Algebra $U\sb {q,p}(\widehat{
\mathfrak{sl}}\sb 2)$: Drinfeld Currents and Vertex Operators}},
   	 {{\it Comm. Math. Phys. }}
  	 {\bf 199}
  	 {(1999)}
  	 {605--647}



\bibitem{Kac}
   {V. G. Kac},
  Infinite Dimensional Lie algebras,
  {\it 3rd. ed. Cambridge University Press},
 1990.



\bibitem{KR}
   {E. Koelink and H. Rosengren},
{{Harmonic Analysis on the  $SU(2)$ Dynamical Quantum Group}},
  {\it Acta.Appl.Math.},
  {\bf 69},
  2001,
  163--220.


\bibitem{KK03}
{T. Kojima and H. Konno},
{{The Elliptic Algebra $U_{q,p}(\slnh)$ and the Drinfeld Realization of the Elliptic 
Quantum Group $\Bqla(\slnh)$}},
{\em Comm.Math.Phys.}{\bf 239} (2003) 405-447. 	 



\bibitem{KKW}
{T. Kojima and H. Konno and R. Weston},
{{The Vertex-Face Ccorrespondence
   and Correlation Functions of the Fusion Eight-Vertex Models I: The General Formalism}},
{\em Nucl.Phys.}{\bf B720} (2005) 348-398.  


\bibitem{K98}
  {H. Konno},
  {An Elliptic Algebra $U_{q,p}(\slth)$ and the Fusion RSOS Models},
  {\it Comm. Math. Phys.}
  {\bf 195}
{(1998)}
   {373--403}.


\bibitem{Konno06}
  {H. Konno},
  {Dynamical $R$ Matrices of Elliptic Quantum Groups and Connection Matrices 
  for the $q$-KZ Equations},
  {\it SIGMA} {\bf 2} (2006) Paper 091.
  
  
\bibitem{Konno08}
  {H. Konno},
  {Elliptic Quantum Group $U_{q,p}(\slth)$ and Vertex Operators},
  {\it J.Phys.A}{\bf 41} (2008) 194012.  	 
  
  
\bibitem{Konno09}
  {H. Konno},
  {Elliptic Quantum Group $U_{q,p}(\slth)$, Hopf Algebroid Structure and Elliptic 
  Hypergoemetric Series},
  {\it J. Geom. Phys.} {\bf 59} (2009) 1485-1511.


\bibitem{Konno1314}
H.Konno, 
\newblock {Elliptic Quantum Group $U_{q,p}(\gh)$ and Deformed $W$-algebras},  talk given at 
{\it Elliptic Integrable Systems and Hypergeometric Functions}, Lorentz Center, 2013; 
%
  {Elliptic Quantum Group, 
Drinfeld Coproduct and Deformed $W$-Algebras}, talk given at 
{\it RAQIS 2014}, Dijion.  


\bibitem{Konno16}
H.Konno,
\newblock An Elliptic Quantum Groups $U_{q,p}(\glnh)$ and $E_{q,p}(\glnh)$,
Preprint 2016, arXiv :1603.04129, to appear in Adv.Studies in Pure Math.  

\bibitem{LB}
H.Lange and Ch.Birkenhake,
\newblock Complex Abelian Varieties, 
Springer-Verlag, 1992.

\bibitem{Matsuo}
A.Matsuo,
\newblock Quantum Algebra Structure of Certain Jackson Integrals,
{\it Comm.Math.Phys.} {\bf 157} (1993) 479--498. 


\bibitem{MO}
D. Maulik and  A. Okounkov,
\newblock Quantum Groups and Quantum Cohomology , Preprint 2012, arXiv:1211.1287, 

\bibitem{Mimachi}
K.Mimachi, 
\newblock A Solution to Quantum Knizhnik-Zamolodchikov Equations and Its Application to Eigenvalue Problems of the Macdonald Type,
{\it Duke Math.J.} {\bf 85} (1996) 635--658. 


\bibitem{MN}
 K.Mimachi and M.Noumi,
\newblock Representations of the Hecke Algebra on a Family of Rational Functions, 
Preprint 1996, unpublished. 


  	 
\bibitem{Nakajima}
H.Nakajima,
\newblock Quiver Varieties and Finite-dimensional Representations of Quantum Affine
 Algebras,
 {\it J. Amer. Math. Soc.}, {\bf 14} (2001) 145--238.


\bibitem{Nakayashiki}
A.Nakayashiki,
\newblock Trace construction of a Basis for the Solution Space of $sl_N$ $q$-KZ Equation.
{\it Comm.Math.Phys.} {\bf 212} (2000) 29--61. 


\bibitem{Negut}
A.Negut,
\newblock The Shuffle Algebra Revisited, 
{\it Int.Math.Res.Not.} {\bf 22} (2014) 6242--6275;  
\newblock Quantum Toroidal and Shuffle Algebras, $R$-matrices and a Conjecture of Kuznetsov,
Preprint 2013, arXiv:1302.6202. 

\bibitem{Okounkov}
A. Okounkov, 
\newblock Lectures on K-theoretic Computations in Enumerative Geometry. (2015),  arXiv: 1512.07363.


\bibitem{RTV}
R.Rim\' anyi, V.Tarasov  and A.Varchenko, 
\newblock Trigonometric Weight Functions as $K$-theoretic Stable Envelope Maps for the Cotangent Bundle of a Flag Variety,
 {\it J.Geom.Phys.},
  {\bf 94} (2015) 81--119.

\bibitem{RV}
R.Rim\' anyi and A.Varchenko, 
\newblock Dynamical Gelfand-Tetlin Algebra and Equivarant Cohomology of Grassmanians, Preprint 2015, 
arXiv:1510.03625. 

\bibitem{RTV17}
R.Rim\' anyi, V.Tarasov  and A.Varchenko, 
\newblock Elliptic and $K$-theoretic Stable Envelopes and Newton Polytopes, Preprint 2017, arXiv:1705.09344. 

\bibitem{TV95}
V.Tarasov and A.Varchenko,
\newblock Jackson Integral Representations for Solutions to the Quantized Knizhnik-Zamolodchikov Equation,
{\it St.Petersburg Math.J.} {\bf 6} (1995) 275-313.  

\bibitem{TV97}
V.Tarasov and A.Varchenko,
\newblock Geometry of $q$-Hypergeometric Functions, Quantum Affine Algebras and Elliptic Quantum Groups,
{\it Ast\' erisque} {\bf 246} (1997) Soci\' et\' e Math\' ematique de France. 



\end{thebibliography}
\end{document}